\documentclass{article}
\usepackage{amssymb}
\usepackage{amsbsy}
\begin{document}
\centerline{\Large \bf } \vskip 6pt

\begin{center}{\Large \bf The dual Orlicz-Aleksandrov-Fenchel
inequality}\end{center}

\vskip 6pt\begin{center} \centerline{Chang-Jian Zhao}
\centerline{\it Department of Mathematics, China Jiliang
University, Hangzhou 310018, P. R. China}\centerline{\it
Email: chjzhao@163.com~~  chjzhao@cjlu.edu.cn}
\end{center}

\vskip 10pt

\begin{center}
\begin{minipage}{12cm}
{\bf Abstract}~ By calculating the first order variation of the dual mixed volumes, we put forward a new concept {\it Orlicz multiple dual mixed volumes.} The classical dual mixed volumes and dual Aleksandrov-Fenchel
inequality are generalized to the Orlicz space. We establish dual Orlicz Aleksandrov-Fenchel inequality and
a new concept of $L_p$ multiple dual mixed volume and $L_p$ dual Aleksandrov-Fenchel inequality are first derived here, too. The dual Aleksandrov-Fenchel inequality and the $L_p$-dual Aleksandrov-Fenchel inequality are all special cases of the dual Orlicz-Aleksandrov-Fenchel inequality.  As an application, a new dual
Orlicz-Brunn-Minkowski inequality for the Orlicz harmonic addition
is also established.


\end{minipage}
\end{center}
\vskip 20pt

\noindent{\large \bf 1 ~Introducation}\vskip 10pt
It is well known that the vector addition is one of the most important operators in geometry.
As an operation between sets $K$ and $L$, defined by
$$K+L=\{x+y: x\in K,y\in L\},$$ it is usually called Minkowski
addition and combine volume play an important role in the
Brunn-Minkowski theory. During the last few decades, the theory
has been extended to $L_{p}$-Brunn-Minkowski theory. The first, a
set called as $L_{p}$ addition, introduced by Firey in [6] and
[7]. Denoted by $+_{p}$, and defined by
$$h(K+_{p}L,x)^{p}=h(K,x)^{p}+h(L,x)^{p},$$
for $p\geq 1$, $x\in {\Bbb R}^{n}$ and compact convex sets $K$ and $L$ in
${\Bbb R}^{n}$ containing the origin. Here the functions are the
support functions. If $K$ is a nonempty closed (not necessarily
bounded) convex set in ${\Bbb R}^{n}$, then
$$h(K,x)=\max\{x\cdot y: y\in K\},$$ for $x\in {\Bbb R}^{n},$ defined the support function $h(K,x)$ of
$K$. A nonempty closed convex set is uniquely determined by its
support function. $L_{p}$ addition and inequalities are the
fundamental and core content in the $L_{p}$ Brunn-Minkowski
theory. For recent important results and more information from
this theory, we refer to x[16] [20-23], [29], [31], [35-39], [42-43],
[48-49], [51-52], x[59] and the references therein.

In recent years, progress towards an Orlicz-Brunn-Minkowski
theory, initiated by Lutwak, Yang and Zhang [40] and [41].
Gardner, Hug and Weil [11] introduced a corresponding addition and
constructed a general framework for the Orlicz-Brunn-Minkowski
theory, and made clear for the first time the relation to Orlicz
spaces and norms, and Orlicz-Minkowski and Brunn-Minkowski
inequalities were established. The Orlicz addition of convex
bodies was also introduced in different ways and extend the
$L_p$-Brunn-Minkowski inequality to the Orlicz-Brunn-Minkowski
inequality (see [53]). The Orlicz centroid inequality for star
bodies was introduced in [61] which is an extension from convex to
star bodies. The other articles advance the theory can be found in
literatures [19], [25], [27], [45] and [55].

The radial addition $K\widetilde{+}L$ of star sets (compact sets
that is star-shaped at $o$ and contains $o$) $K$ and $L$ can be
defined by
$$\rho(K\widetilde{+}L,\cdot)=\rho(K,\cdot)+\rho(L,\cdot),$$
where $\rho(K,\cdot)$ denotes the radial function of star set $K.$ The
radial function of star set $K$ is defined by
$$\rho(K,u)=\max\{c\geq 0: cu\in K\},$$ for $u\in S^{n-1}$.
The origin and history of the radial addition can be referred to [10, p. 235].
When $\rho(K,\cdot)$ is positive and continuous, $K$
will be called a star body. Let ${\cal S}^{n}$ denote the set of
star bodies about the origin in ${\Bbb R}^{n}$. The radial addition and volume are the core and essence of the classical dual Brunn-Minkowski theory and played an important role in theory (see, e.g., [3], [14], [17], [18], [26], [28], [32] and [51] for recent important contributions). Lutwak [34] introduced the concept of dual
mixed volume laid the foundation of the dual Bruun-Minkowski
theory. What is particularly important is that this theory plays a very important and key role in solving the Busemann-Petty problem in
[8], [14], [33] and [47]. For any $p\neq
0$, the $L_p$-radial addition $K\widetilde{+}_{p}L$ defined by
(see [9])
$$\rho(K\widetilde{+}_{p}L,x)^{p}=\rho(K,x)^{p}+\rho(L,x)^{p},$$
for $x\in {\Bbb R}^{n}$ and $K,L\in{\cal S}^{n}$ (see [12]). Obviously, when $p=1$, the $L_p$-radial addition $\widetilde{+}_{p}$ becomes the radial addition $\widetilde{+}$.

The $L_p$-harmonic radial addition was defined by Lutwak [31]: If $K,L$ are star bodies,
the $L_p$-harmonic radial addition, defined by
$$\rho(K\widehat{+}_{p}L,x)^{-p}=\rho(K,x)^{-p}+\rho(L,x)^{-p},\eqno(1.1)$$
for $p\geq 1$ and $x\in {\Bbb R}^{n}$. The $L_p$-harmonic radial addition of convex bodies was first studied by
Firey [6]. The operation of the $L_p$-harmonic radial addition and the
$L_{p}$-dual Minkowski, Brunn-Minkwski inequalities are the basic concept and inequalities in the $L_{p}$-dual
Brunn-Minkowski theory. The latest information and important results of this theory can be referred to [3], [5], [9], [12], [15], [17-18] and [28] and the references therein.
A systematic investigation on the concepts of the addition for convex body and star body, we refer to [11], [12] and [15]. The theory has recently turned to a study
extending from $L_p$-dual Brunn-Minkowski theory to the Orlicz
dual Brunn-Minkowski theory. The dual Orlicz-Brunn-Minkowski
theory has also attracted the attention, see [13], [54], [56-58] and [60]. The Orlicz harmonic radial addition
$K\widehat{+}_{\phi}L$ of two star bodies $K$ and $L$, defined by (see [60])
$$\rho(K\widehat{+}_{\phi}L,u))=\sup\left\{\lambda>0:\phi\left(\frac{\rho(K,u)}{\lambda}\right)+\phi\left(\frac{\rho(L,u)}{\lambda}\right)
\leq\phi(1)\right\},\eqno(1.2)$$ where $u\in
S^{n-1}$, and $\phi:(0,\infty)\rightarrow(0,\infty)$ is a convex and decreasing
function such that $\phi(0)=\infty$,
$\lim_{t\rightarrow\infty}\phi(t)=0$ and
$\lim_{t\rightarrow0}\phi(t)=\infty.$ let ${\cal C}$ denote the
class of the convex and decreasing functions $\phi$. Obviously, if
${\phi}(t)=t^{-p}$ and $p\geq 1$, then the Orlicz harmonic addition
becomes the $L_p$-harmonic radial addition. The Orlicz dual mixed volume, denoted by
$\widetilde{V}_{\phi}(K,L)$, defined by
$$\widetilde{V}_{\phi}(K,L):=\frac{{\phi}'_{r}(1)}{n}\lim_{\varepsilon\rightarrow
0^+}\frac{V(K\widehat{+}_{\phi}\varepsilon\cdot
L)-V(K)}{\varepsilon}=\frac{1}{n}\int_{S^{n-1}}\phi\left(\frac{\rho(L,u)}{\rho(K,u)}\right)
\rho(K,u)^{n}dS(u),\eqno(1.3)$$ where
$K\widehat{+}_{\phi}\varepsilon\cdot L$ is the Orlicz linear
combination of $K$ and $L$ (see Section 3), and the right
derivative of a real-valued function $\phi$ is denoted by $\phi'_{r}$.
When ${\phi}(t)=t^{-p}$ and $p\geq 1$, the Orlicz dual mixed volume
$\widetilde{V}_{\phi}(K,L)$ becomes the $L_p$ dual mixed volume
$\widetilde{V}_{-p}(K,L)$, defined by (see [31])
$$\widetilde{V}_{-p}(K,L)=\frac{1}{n}\int_{S^{n-1}}\rho(K,u)^{n+p}\rho(L,u)^{-p}dS(u).\eqno(1.4)$$

If $K_{1},\ldots,K_{n}\in {\cal S}^{n}$, the dual mixed volume of
star bodies $K_{1},\ldots,K_{n},$
$\widetilde{V}(K_{1},\ldots,K_{n})$ defined by Lutwak (see [34])
$$\widetilde{V}(K_{1},\ldots,K_{n})=\frac{1}{n}\int_{S^{n-1}}\rho(K_{1},u)\cdots\rho(K_{n},u)dS(u).\eqno(1.5)$$
Lutwak's dual Aleksandrov-Fenchel inequality (see [33]) is
following: If $K_{1},\cdots,K_{n}\in {\cal S}^{n}$ and $1\leq
r\leq n$, then
$$\widetilde{V}(K_{1},\cdots,
K_{n})\leq\prod_{i=1}^{r}\widetilde{V}(K_{i}\ldots,K_{i},K_{r+1},\ldots,K_{n})^{\frac{1}{r}},$$
with equality if and only if $K_{1},\ldots,K_{r}$ are all
dilations of each other.

As we all know, $L_{p}$-dual volume $\widetilde{V}_{p}(K,L)$ of
star bodies $K$ and $L$ has been extended to the Orlicz space and
becomes Orlicz dual volume $\widetilde{V}_{\phi}(K,L)$. However
the dual mixed volumes $\widetilde{V}(K_{1},\cdots,K_{n})$ has not
been extended to the Orlicz space, and this question becomes a
difficult research in convex geometry. The main reason for the
difficulty is that it has not been extended to $L_p$ space. In the
paper, we solved the difficulty. Our main aim is to generalize
direct the classical dual mixed volumes
$\widetilde{V}(K_{1},\cdots,K_{n})$ and dual Aleksandrov-Fenchel
inequality to the Orlicz space without passing through $L_{p}$
space. Under the framework of dual Orlicz-Brunn-Minkowski theory,
we introduce a new affine geometric quantity by calculating the
first order Orlicz variation of the dual mixed volumes, and call
it Orlicz multiple dual mixed volumes, denoted by
$\widetilde{V}_{\phi}(K_{1},\cdots,K_{n},L_{n})$, which involves
$(n+1)$ star bodies in ${\Bbb R}^{n}$. The fundamental notions and
conclusions of the dual mixed volume
$\widetilde{V}(K_{1},\cdots,K_{n})$ and the dual Minkowski, and
Aleksandrov-Fenchel inequalities are extended to an Orlicz
setting. The related concepts and conclusions of $L_p$ multiple
dual mixed volume $\widetilde{V}_{p}(K_{1},\cdots,K_{n},L_{n})$
and $L_p$-Aleksandrov-Fenchel inequality are also derived. The new
dual Orlicz-Aleksandrov-Fenchel inequality in special case yields
the dual Aleksandrov-Fenchel inequality and the Orlicz dual
Minkowski inequality for the dual quermassintegrals, respectively.
As application, a new dual Orlicz-Brunn-Minkowski inequality for
the Orlicz harmonic radial addition is established, which implies
the dual Orlicz-Brunn-Minkowski inequality for the dual
quermassintegrals.

Complying with the spirit of introduction of Aleksandrov, Fenchel
and Jensen's mixed quermassintegrals, and introduction of Lutwak's
$L_p$-mixed quermassintegrals,  we calculate the first order
Orlicz variational of the dual mixed volumes. In Section 4, we
prove that the first order Orlicz variation of the dual mixed
volumes can be expressed as:
$$\frac{d}{d\varepsilon}\bigg|_{\varepsilon=0^{+}}\widetilde{V}(L_{1}\widehat{+}_{\phi}\varepsilon\cdot
K_{1},K_{2},\cdots,K_{n})=\frac{1}{{\phi}'_{r}(1)}\cdot\widetilde{V}_{{\phi}}(L_{1},K_{1},\cdots,K_{n}),$$
where $L_{1},K_{1},\cdots,K_{n}\in {\cal S}^{n}$, ${\phi}\in {\cal
C}$ and $\varepsilon>0$. In the above first order variational
equation, we find a new geometric quantity. Based on this, we
extract the required geometric quantity, denoted by
$\widetilde{V}_{{\phi}}(L_{1},K_{1},\cdots,K_{n})$ and call it
Orlicz multiple dual mixed volume
$\widetilde{V}_{{\phi}}(L_{1},K_{1},\cdots,K_{n})$ of $(n+1)$ star
bodies $L_{1},K_{1},\cdots,K_{n}$, defined by
$$\widetilde{V}_{{\phi}}(L_{1},K_{1},\cdots,K_{n}):={\phi}'_{r}(1)\cdot\frac{d}{d\varepsilon}\bigg|_{\varepsilon=0^{+}}\widetilde{V}(L_{1}\widehat{+}_{\phi}\varepsilon\cdot
K_{1},K_{2},\cdots,K_{n}).$$ We also prove the new affine
geometric quantity
$\widetilde{V}_{{\phi}}(L_{1},K_{1},\cdots,K_{n})$ has an integral
representation.
$$\widetilde{V}_{{\phi}}(L_{1},K_{1},\cdots,K_{n})=\frac{1}{n}\int_{S^{n-1}}{\phi}
\left(\frac{\rho(K_{1},u)}{\rho(L_{1},u)}\right)\rho(L_{1},u)\rho(K_{2},u)\cdots\rho(K_{n},u)
dS(u),~~~~~~~~\eqno(1.6)$$ When $K_{1}=L$ and
$L_{1}=K_{2}=\cdots=K_{n}=K$, the Orlicz multiple dual mixed
volume $\widetilde{V}_{{\phi}}(L_{1},K_{1},\cdots,K_{n})$ becomes
the usual Orlicz dual mixed volume $\widetilde{V}_{\phi}(K,L).$
When $\phi(t)=t^{-p},$ $p\geq 1$, Orlicz multiple dual mixed
volume $\widetilde{V}_{{\phi}}(L_{1},K_{1},\cdots,K_{n})$ becomes
a new dual mixed volume in $L_{p}$ place, denoted by
$\widetilde{V}_{-p}(K_{1},\ldots,K_{n},L_{n})$, call it $L_{p}$
multiple dual mixed volume. From (1.6), we have
$$\widetilde{V}_{-p}(L_{1},K_{1},\cdots,K_{n})=\frac{1}{n}\int_{S^{n-1}}\rho(K_{1},u)^{-p}\rho(L_{1},u)^{1+p}\rho(K_{2},u)\cdots\rho(K_{n},u)
dS(u).\eqno(1.7)$$ Putting $L_{1}=K_{1}$ in (1.7), the $L_{p}$
multiple dual mixed volume
$\widetilde{V}_{-p}(L_{1},K_{1},\cdots,K_{n})$ becomes the usual
dual mixed volume $\widetilde{V}(K_{1},\cdots,K_{n}).$ Putting
$K_{1}=L$ and $L_{1}=K_{2}=\cdots=K_{n}=K$ in (1.7),
$\widetilde{V}_{-p}(L_{1},K_{1},\cdots,K_{n})$ becomes the $L_p$
dual mixed volume $\widetilde{V}_{-p}(K,L).$ Putting $K_{1}=L$ and
$L_{1}=K_{2}=\cdots=K_{n-i}=K$ and $K_{n-i+1}=\cdots=K_{n}=B$ in
(1.7), $\widetilde{V}_{-p}(L_{1},K_{1},\cdots,K_{n})$ becomes the
harmonic mixed $p$-quermassintegral, $\widetilde{W}_{-p,i}(K,L)$,
defined by (see Section 2)
$$\widetilde{W}_{-p,i}(K,L)=\frac{1}{n}\int_{S^{n-1}}\rho(K,u)^{n-i+p}\rho(L,u)^{-p}dS(u).\eqno(1.8)$$

In Section 5, we establish the following dual
Orlicz-Aleksandrov-Fenchel inequality for the Orlicz multiple dual
mixed volume. If $L_{1},K_{1},\cdots,K_{n}\in {\cal S}^{n}$,
${\phi}\in {\cal C}$ and $1\leq r\leq n$, then
$$\widetilde{V}_{\phi}(L_{1},K_{1},K_{2},\cdots,
K_{n})\geq\widetilde{V}(L_{1},K_{2},\cdots, K_{n})\cdot
\phi\left(\frac{\prod_{i=1}^{r}\widetilde{V}(K_{i}\ldots,K_{i},K_{r+1},\ldots,K_{n})^{\frac{1}{r}}}
{\widetilde{V}(L_{1},K_{2}\ldots,K_{n})}\right).\eqno(1.9)$$ If
$\phi$ is strictly convex, equality holds if and only if
$L_{1},K_{1},\ldots,K_{r}$ are all dilations of each other. When
$\phi(t)=t^{-p}$, $p=1$ and $K_{1}=L_{1}$, the dual
Orlicz-Aleksandrov-Fenchel inequality becomes Lutwak's dual
Aleksandrov-Fenchel inequality. If $K_{1},\cdots,K_{n}\in {\cal
S}^{n}$ and $1\leq r\leq n$, then
$$\widetilde{V}(K_{1},\cdots,
K_{n})\leq\prod_{i=1}^{r}\widetilde{V}(K_{i}\ldots,K_{i},K_{r+1},\ldots,K_{n})^{\frac{1}{r}},\eqno(1.10)$$
with equality if and only if $K_{1},\ldots,K_{r}$ are all
dilations of each other. When $\phi(t)=t^{-p}$, $p\geq 1$, the
dual Orlicz-Aleksandrov-Fenchel inequality (1.9) becomes the
following $L_p$ dual Aleksandrov-Fenchel inequality. If
$L_{1},K_{1},\cdots,K_{n}\in {\cal S}^{n}$, $p\geq 1$ and $1\leq
r\leq n$, then
$$\frac{\widetilde{V}(L_{1},K_{2},\cdots,
K_{n})^{p+1}}{\widetilde{V}_{-p}(L_{1},K_{1},K_{2},\cdots,
K_{n})}\leq\prod_{i=1}^{r}\widetilde{V}(K_{i}\ldots,K_{i},K_{r+1},\ldots,K_{n})^{\frac{p}{r}}.\eqno(1.11)$$
If $\phi$ is strictly convex, equality holds if and only if
$L_{1},K_{1},\ldots,K_{r}$ are all dilations of each other. When
$K_{1}=L$ and $L_{1}=K_{2}=\cdots=K_{n}=K$, the Orlicz dual
Aleksandrov-Fenchel inequality (1.9) becomes the following Orlicz
dual Minkowski inequality established in [60]. If $K,L\in {\cal
S}^{n}$ and ${\phi}\in{\cal C}$, then
$$\widetilde{V}_{\phi}(K,L)\geq V(K)\cdot
\phi\left(\left(\frac{V(L)}{V(K)}\right)^{\frac{1}{n}}\right),\eqno(1.12)$$
with equality if and only if $K$ and $L$ are dilates. In fact, in
Section 5, we show also the Orlicz Aleksandrov-Fenchel inequality
(1.9) in special case yield also the following result. If $K,L\in
{\cal S}^{n}$, $0\leq i<n$ and $\phi\in {\cal C}$, then
$$\widetilde{W}_{\phi,i}(K,L)\geq
\widetilde{W}_{i}(K)\phi\left(\left(\frac{\widetilde{W}_{i}(L)}{\widetilde{W}_{i}(K)}\right)^{1/(n-i)}\right)
.\eqno(1.13)$$ If $\phi$ is strictly convex, equality holds if and
only if $K$ and $L$ are dilates. Here $\widetilde{W}_{i}(K)$ is
the usually dual quermassintegral of $K$, and
$\widetilde{W}_{\phi,i}(K,L)$ is the Orlicz dual mixed
quermassintegral of $K$ and $L$, defined by (see Section 4)
$$\widetilde{W}_{\phi,i}(K,L)=\frac{1}{n}\int_{S^{n-1}}\phi
\left(\frac{\rho(L,u)}{\rho(K,u)}\right)\rho(K,u)^{n-i}dS(u).\eqno(1.14)$$

In Section 6, we establish the following Orlicz dual
Brunn-Minkowski type inequality. If $L_{1},K_{1},\cdots,K_{n}\in
{\cal S}^{n}$ and ${\phi}\in {\cal C}$, then
$$\phi(1)\geq\phi\left(\left(\frac{V(K_{1})\cdots
V(K_{n})}{\widetilde{V}(K_{1}\widehat{+}_{\phi}L_{1},K_{2},\cdots,
K_{n})^{n}}\right)^{\frac{1}{n}}\right)+\phi\left(\left(\frac{V(L_{1})V(K_{2})\cdots
V(K_{n})}{\widetilde{V}(K_{1}\widehat{+}_{\phi}L_{1},K_{2},\cdots,
K_{n})^{n}}\right)^{\frac{1}{n}}\right).\eqno(1.15)$$ If $\phi$ is
strictly convex, equality holds if and only if
$L_{1},K_{1},\ldots,K_{n}$ are all dilations of each other.
Putting $K_{1}=K, L_{1}=L$ and
$K_{2}=\cdots=K_{n}=K_{1}\widehat{+}_{\phi}L_{1}$ in (1.15), it
follows the Orlicz dual Brunn-Minkowski inequality established in
[58]. If $K,L\in {\cal S}^{n}$ and ${\phi}\in{\cal C}$
$$\phi(1)\geq\phi\left(\left(\frac{V(K)}
{V(K\widehat{+}_{\phi}L)}\right)^{\frac{1}{n}}\right)+
\phi\left(\left(\frac{V(L)}{V(K\widehat{+}_
{\phi}L)}\right)^{\frac{1}{n}}\right).\eqno(1.16)$$ If $\phi$ is
strictly convex, equality holds if and only if $K$ and $L$ are
dilates. In fact, in Section 6, we show also the Orlicz dual
Brunn-Minkowski inequality (1.15) in special case yield also the
following result. If $K,L\in {\cal S}^{n}$, ${\phi}\in {\cal C}$
and $0\leq i<n-1$, then
$$\phi(1)\geq\phi\left(\left(\frac{\widetilde{W}_{i}(K)}{\widetilde{W}_{i}(K\widehat{+}_{\phi}L)}\right)^{\frac{1}{n-i}}\right)+\phi\left(\left(\frac{\widetilde{W}_{i}(L)}
{\widetilde{W}_{i}(K\widehat{+}_{\phi}L)}\right)^{\frac{1}{n-i}}\right).\eqno(1.17)$$
If $\phi$ is strictly convex, equality holds if and only if $K$
and $L$ are dilates.

\vskip 10pt \noindent{\large \bf 2 ~Preliminaries}\vskip 10pt

The setting for this paper is $n$-dimensional Euclidean space
${\Bbb R}^{n}$. A body in ${\Bbb R}^{n}$ is a compact set equal to
the closure of its interior. For a compact set $K\subset {\Bbb
R}^{n}$, we write $V(K)$ for the ($n$-dimensional) Lebesgue
measure of $K$ and call this the volume of $K$. The unit ball in
${\Bbb R}^{n}$ and its surface are denoted by $B$ and $S^{n-1}$,
respectively. Let ${\cal K}^{n}$ denote the class of nonempty
compact convex subsets containing the origin in their interiors in
${\Bbb R}^{n}$. Associated with a compact subset $K$ of ${\Bbb
R}^n$, which is star-shaped with respect to the origin and
contains the origin, its radial function is $\rho(K,\cdot):
S^{n-1}\rightarrow [0,\infty),$ defined by
$$\rho(K,u)=\max\{\lambda\geq 0: \lambda u\in K\}.$$
Two star bodies $K$ and $L$ are dilates if $\rho(K,u)/\rho(L,u)$
is independent of $u\in S^{n-1}$. If $\lambda>0$, then
$$\rho(\lambda K,u)=\lambda\rho(K,u).$$
From the definition of the radial function, it follows immediately
that for $A\in GL(n)$ the radial function of the image
$AK=\{Ay:y\in K\}$ of $K$ is given by (see e.g. [10])
$$\rho(AK,u)=\rho(K, A^{-1}u),$$
for all $u\in S^{n-1}$. Namely, the radial function is homogeneous
of degree $-1$. Let $\tilde{\delta}$ denote the radial Hausdorff
metric, as follows, if $K, L\in {\cal S}^{n}$, then (see e.g.
[46])
$$\tilde{\delta}(K,L)=|\rho(K,u)-\rho(L,u)|_{\infty}.$$

\vskip 8pt {\it 2.1~ Dual mixed volumes}\vskip 8pt

The polar coordinate formula for volume of a compact set $K$ is
$$V(K)=\frac{1}{n}\int_{S^{n-1}}\rho(K,u)^{n}dS(u).\eqno(2.1)$$ The
first dual mixed volume, $\widetilde{V}_{1}(K,L)$, defined by
$$\widetilde{V}_{1}(K,L)=\frac{1}{n}\lim_{\varepsilon\rightarrow
0^{+}}\frac{V(K\widetilde{+}\varepsilon\cdot
L)-V(K)}{\varepsilon},$$ where $K,L\in{\cal S}^{n}.$ The integral
representation for first dual mixed volume is proved: For
$K,L\in{\cal S}^{n},$
$$\widetilde{V}_{1}(K,L)=\frac{1}{n}\int_{S^{n-1}}\rho(K,u)^{n-1}\rho(L,u)dS(u).\eqno(2.2)$$
The Minkowski inequality for first dual mixed volume is the
following: If $K,L\in{\cal S}^{n},$ then
$$\widetilde{V}_{1}(K,L)^{n}\leq V(K)^{n-1}V(L),$$
with equality if and only if $K$ and $L$ are dilates. (see [33])
If $K_{1},\ldots,K_{n}\in {\cal S}^{n}$, the dual mixed volume
$\widetilde{V}(K_{1},\ldots,K_{n})$ is defined by (see [34])
$$\widetilde{V}(K_{1},\ldots,K_{n})=\frac{1}{n}\int_{S^{n-1}}\rho(K_{1},u)\cdots\rho(K_{n},u)dS(u).\eqno(2.3)$$
If $K_{1}=\cdots=K_{n-i}=K,$ $K_{n-i+1}=\cdots=K_{n}=L$, the dual
mixed volume $\widetilde{V}(K_{1},\ldots,K_{n})$ is written as
$\widetilde{V}_{i}(K,L)$. If $L=B,$ the dual mixed volume
$\widetilde{V}_{i}(K,L)=\widetilde{V}_{i}(K,B)$ is written as
$\widetilde{W}_{i}(K)$ and called dual quermassintegral of $K$.
For $K\in {\cal S}^{n}$ and $0\leq i<n$,
$$\widetilde{W}_{i}(K)=\frac{1}{n}\int_{S^{n-1}}\rho(K,u)^{n-i}dS(u).\eqno(2.4)$$
If $K_{1}=\cdots=K_{n-i-1}=K,$ $K_{n-i}=\cdots=K_{n-1}=B$ and
$K_{n}=L$, the dual mixed volume
$\widetilde{V}(\underbrace{K,\ldots,K}_{n-i-1},\underbrace{B,\ldots,B}_{i},L)$
is written as $\widetilde{W}_{i}(K,L)$ and called dual mixed
quermassintegral of $K$ and $L$. For $K,L\in {\cal S}^{n}$ and
$0\leq i<n$, it is easy that ([45])
$$\widetilde{W}_{i}(K,L)=\lim_{\varepsilon\rightarrow 0^{+}}\frac{\widetilde{W}_{i}(K\widetilde{+}
\varepsilon\cdot
L)-\widetilde{W}_{i}(K)}{\varepsilon}=\frac{1}{n}\int_{S^{n-1}}\rho(K,u)^{n-i-1}\rho(L,u)dS(u).\eqno(2.5)$$

The fundamental inequality for dual mixed quermassintegral stated
that: If $K,L\in {\cal S}^{n}$ and $0\leq i<n$, then
$$\widetilde{W}_{i}(K,L)^{n-i}\leq\widetilde{W}_{i}(K)^{n-1-i}\widetilde{W}_{i}(L),\eqno(2.6)$$
with equality if and only if $K$ and $L$ are dilates. The
Brunn-Minkowski inequality for dual quermassintegral is the
following: If $K,L\in {\cal S}^{n}$ and $0\leq i<n$, then
$$\widetilde{W}_{i}(K\widetilde{+}L)^{1/(n-i)}\leq\widetilde{W}_{i}(K)^{1/(n-i)}+\widetilde{W}_{i}(L)^{1/(n-i)},\eqno(2.7)$$
with equality if and only if $K$ and $L$ are dilates.

\vskip 8pt {\it 2.2~ $L_{p}$-dual mixed volume}\vskip 8pt

The dual mixed volume $\widetilde{V}_{-1}(K,L)$ of star bodies $K$
and $L$ is defined by ([31])
$$\widetilde{V}_{-1}(K,L)=\lim_{\varepsilon\rightarrow
0^{+}}\frac{V(K)-V(K\widehat{+}\varepsilon\cdot
L)}{\varepsilon},\eqno(2.8)$$ where $\widehat{+}$ is the harmonic
addition. The following is a integral representation for the dual
mixed volume $\widetilde{V}_{-1}(K,L)$:
$$\widetilde{V}_{-1}(K,L)=\frac{1}{n}\int_{S^{n-1}}\rho(K,u)^{n+1}\rho(L,u)^{-1}dS(u).\eqno(2.9)$$
The dual Minkowski inequality for the dual mixed volume states
that
$$\widetilde{V}_{-1}(K,L)^{n}\geq V(K)^{n+1} V(L)^{-1},\eqno(2.10)$$
with equality if and only if $K$ and $L$ are dilates. (see ([32]))

The dual Brunn-Minkowski inequality for the harmonic addition
states that
$$V(K\widehat{+}L)^{-1/n}\geq V(K)^{-1/n}+V(L)^{-1/n},\eqno(2.11)$$
with equality if and only if $K$ and $L$ are dilates (This
inequality is due to Firey [6]).

The $L_p$ dual mixed volume $\widetilde{V}_{-p}(K,L)$ of $K$ and
$L$ is defined by ([31])
$$\widetilde{V}_{-p}(K,L)=-\frac{p}{n}\lim_{\varepsilon\rightarrow 0^{+}}\frac{V(K\widehat{+}_{p}\varepsilon\cdot L)-V(K)}{\varepsilon},\eqno(2.12)$$
where $K,L\in {\cal S}^{n}$ and $p\geq 1$.

The following is an integral representation for the $L_p$ dual
mixed volume: For $K,L\in {\cal S}^{n}$ and $p\geq 1$,
$$\widetilde{V}_{-p}(K,L)=\frac{1}{n}\int_{S^{n-1}}\rho(K,u)^{n+p}\rho(L,u)^{-p}dS(u).\eqno(2.13)$$
$L_{p}$-dual Minkowski and Brunn-Minkowski inequalities were
established by Lutwak [31]: If $K,L\in {\cal S}^{n}$ and $p\geq
1$, then
$$\widetilde{V}_{-p}(K,L)^{n}\geq V(K)^{n+p} V(L)^{-p},\eqno(2.14)$$
with equality if and only if $K$ and $L$ are dilates, and
$$V(K\widehat{+}_{p}L)^{-p/n}\geq V(K)^{-p/n}+V(L)^{-p/n},\eqno(2.15)$$
with equality if and only if $K$ and $L$ are dilates.

\vskip 8pt {\it 2.3~ Mixed $p$-harmonic quermassintegral}\vskip
8pt

From (1.1), it is easy to see that if $K,L\in {\cal S}^{n}$,
$0\leq i<n$ and $p\geq 1$, then
$$-\frac{p}{n-i}\lim_{\varepsilon\rightarrow 0^{+}}\frac{\widetilde{W}_{i}(K\widehat{+}_{p}\varepsilon\cdot L)
-\widetilde{W}_{i}(L)}{\varepsilon}=\frac{1}{n}\int_{S^{n-1}}\rho(K.u)^{n-i+p}\rho(L.u)^{-p}dS(u).\eqno(2.16)$$
Let $K,L\in {\cal S}^{n}$, $0\leq i<n$ and $p\geq 1$, the mixed
$p$-harmonic quermassintegral of star $K$ and $L$, denoted by
$\widetilde{W}_{-p,i}(K,L)$, defined by (see [50])
$$\widetilde{W}_{-p,i}(K,L)=\frac{1}{n}\int_{S^{n-1}}\rho(K,u)^{n-i+p}\rho(L,u)^{-p}dS(u).\eqno(2.17)$$
Obviously, when $K=L$, the $p$-harmonic quermassintegral
$\widetilde{W}_{-p,i}(K,L)$ becomes the dual quermassintegral
$\widetilde{W}_{i}(K)$. The Minkowski and Brunn-Minkowski
inequalities for the mixed $p$-harmonic quermassintegral are
following (see [50]): If $K,L\in {\cal S}^{n}$, $0\leq i<n$ and
$p\geq 1$, then
$$\widetilde{W}_{-p,i}(K,L)^{n-i}\geq\widetilde{W}_{i}(K)^{n-i+p}\widetilde{W}_{i}(L)^{-p},\eqno(2.18)$$
with equality if and only if $K$ and $L$ are dilates. If $K,L\in
{\cal S}^{n}$, $0\leq i<n$ and $p\geq 1$, then
$$\widetilde{W}_{i}(K\widehat{+}_{p}L)^{-p/(n-i)}\geq\widetilde{W}_{i}(K)^{-p/(n-i)}+\widetilde{W}_{i}(L)^{-p/(n-i)},\eqno(2.19)$$
with equality if and only if $K$ and $L$ are dilates.

Inequality (2.19) is a Brunn-Minkowski type inequality for the
$p$-harmonic addition. For different variants of the classical
Brunn-Minkowski inequalities we refer to [1], [2], [4], [43] and
[47] and the references therein.

\vskip 10pt \noindent{\large \bf 3 ~Orlicz harmonic linear
combination}\vskip 10pt

Throughout the paper, the standard orthonormal basis for ${\Bbb
R}^{n}$ will be $\{e_{1},\ldots,e_{n}\}$. Let ${\cal C}_{m},$
$m\in{\Bbb N}$, denote the set of convex function
$\phi:[0,\infty)^{m} \rightarrow (0,\infty)$ that are strictly
decreasing in each variable and satisfy $\phi(0)=\infty$. When
$m=1$, we shall write ${\cal C}$ instead of ${\cal C}_{1}$. Orlicz
harmonic radial addition is defined below.

{\bf Definition 3.1}~ Let $m\geq 2,\phi\in{\cal C}_{m}$, $K_{j}\in
{\cal S}^{n}$ and $j=1,\ldots,m$, define the Orlicz harmonic
addition of $K_{1},\ldots,K_{m}$, denoted by
$\widehat{+}_{\phi}(K_{1},\ldots,K_{m})$, defined by
$$\rho(\widehat{+}_{\phi}(K_{1},\ldots,K_{m}),x)=\sup\left\{\lambda>0:
\phi\left(\frac{\rho(K_{1},x)}{\lambda},\ldots,\frac{\rho(K_{m},x)}{\lambda}\right)\leq
\phi(1)\right\},\eqno(3.1)$$ for $x\in{\Bbb R}^{n}.$ Equivalently,
the Orlicz multiple harmonic addition
$\widehat{+}_{\phi}(K_{1},\ldots,K_{m})$ can be defined implicitly
by
$$\phi\left(\frac{\rho(K_{1},x)}{\rho(\widehat{+}_{\phi}(K_{1},
\ldots,K_{m}),x)},\ldots,\frac{\rho(K_{m},x)}{\rho(\widehat{+}_{\phi}(K_{1},
\ldots,K_{m}),x)}\right)=\phi(1),\eqno(3.2)$$ for all $x\in {\Bbb
R}^{n}$. An important special case is obtained when
$$\phi(x_{1},\ldots,x_{m})=\sum_{j=1}^{m}\phi(x_{j}),$$
for $\phi(t)\in {\cal C}_{m}$. We then write
$\widehat{+}_{\phi}(K_{1},\ldots,K_{m})=K_{1}\widehat{+}_{\phi}\cdots\widehat{+}_{\phi}K_{m}.$
This means that
$K_{1}\widehat{+}_{\phi}\cdots\widehat{+}_{\phi}K_{m}$ is defined
either by
$$\rho(K_{1}\widehat{+}_{\phi}\cdots\widehat{+}_{\phi}K_{m},x)=
\sup\left\{\lambda>0:\sum_{j=1}^{m}\phi\left(\frac{\rho(K_{j},x)}{\lambda}
\right)\leq \phi(1)\right\},\eqno(3.3)$$ for all $x\in {\Bbb
R}^{n}$, or by the corresponding special case of (3.1). From
(3.3), it follows easy that
$$\sum_{j=1}^{m}\phi\left(\frac{\rho(K_{j},x)}{\lambda}
\right)=\phi(1),$$ if and only if
$$\lambda=\rho(K_{1}\widehat{+}_{\phi}\cdots\widehat{+}_{\phi}K_{m},x).\eqno(3.4)$$

Next, define the Orlicz harmonic linear combination on the case
$m=2$.

{\bf Definition 3.2}~ The Orlicz harmonic linear combination,
denotes $\widehat{+}_{\phi}(K,L,\alpha,\beta)$, defined by
$$\alpha\phi\left(\frac{\rho(K,x)}{\rho(\widehat{+}_{\phi}(K,L,\alpha,\beta),x)}\right)+
\beta\phi\left(\frac{\rho(L,x)}
{\rho(\widehat{+}_{\phi}(K,L,\alpha,\beta),x)}\right)=\phi(1),\eqno(3.5)$$
for $K,L\in {\cal S}^{n}$, $x\in {\Bbb R}^{n}$ and
$\alpha,\beta\geq 0$ (not both zero).

When $\phi(t)=t^{-p}$ and $p\geq 1$, then Orlicz harmonic linear
combination $\widehat{+}_{\phi}(K,L,\alpha,\beta)$ changes to the
$L_p$-harmonic linear combination $\alpha\cdot
K\widehat{+}_{p}\beta\cdot L.$ We shall write
$K\widehat{+}_{\phi}\varepsilon\cdot L$ instead of
$\widehat{+}_{\phi}(K,L,1,\varepsilon)$, for $\varepsilon\geq 0$
and assume throughout that this is defined by (3.5), where
$\alpha=1, \beta=\varepsilon$ and $\phi\in {\cal C}$. It is easy
that $\widehat{+}_{\phi}(K,L,1,1)=K\widehat{+}_{\phi}L.$

\vskip 10pt \noindent{\large \bf 4 ~ Orlicz multiple dual mixed
volume}\vskip 10pt

Let us introduce the Orlicz multiple dual mixed volume.

{\bf Definition 4.1}~ For $\phi\in {\cal C}$, we define the Orlicz
dual mixed volume $\widetilde{V}_{\phi}(L_{1},K_{1},\cdots,K_{n})$
by
$$\widetilde{V}_{\phi}(L_{1},K_{1},\cdots ,K_{n})=:\frac{1}{n}\int_{S^{n-1}}{\phi}
\left(\frac{\rho(K_{1},u)}{\rho(L_{1},u)}\right)\rho(L_{1},u)\rho(K_{2},u)\cdots\rho(K_{n},u)
dS(u),$$ for all $L_{1},K_{1},\ldots,K_{n}\in {\cal S}^{n}$.

To derive this definition 4.1, we need the following lemmas.

{\bf Lemma 4.2}~ ([60]) {\it If $K_{1},L_{1}\in {\cal S}^{n}$ and
$\phi\in {\cal C}$, then}
$$L_{1}\widehat{+}_{\phi}\varepsilon\cdot K_{1}\rightarrow L_{1}\eqno(4.1)$$
as $\varepsilon\rightarrow 0^{+}.$

{\bf Lemma 4.3}~ {\it If $L_{1},K_{1},\ldots,K_{n}\in {\cal
S}^{n}$ and ${\phi}\in {\cal C}$, then}
$$\frac{d}{d\varepsilon}\bigg|_{\varepsilon=0^{+}}\widetilde{V}(L_{1}\widehat{+}_{\phi}\varepsilon\cdot
K_{1},K_{2},\cdots,K_{n})=\frac{1}{n{\phi}'_{r}(1)}\int_{S^{n-1}}{\phi}
\left(\frac{\rho(K_{1},u)}{\rho(L_{1},u)}\right)\rho(L_{1},u)\rho(K_{2},u)\cdots\rho(K_{n},u)
dS(u).\eqno(4.2)$$

{\bf Proof}~ Suppose $\varepsilon>0$, $K_{1},L_{1}\in {\cal
S}^{n}$ and $u\in S^{n-1}$, let
$$\rho_{\varepsilon}=\rho(L_{1}\widehat{+}_{\phi}\varepsilon\cdot K_{1},u).$$
Since
$$\frac{\rho(L_{1},u)}{\rho_{\varepsilon}}=\phi^{-1}\left(\phi(1)-\varepsilon{\phi}
\left(\frac{\rho(K_{1},u)}{\rho_{\varepsilon}}\right)\right),$$
and from Lemma 4.2, and noting that ${\phi}$ is continuous
function, we obtain
$$\lim_{\varepsilon\rightarrow 0^+}\frac{(\rho_{\varepsilon}-\rho(L_{1},u))\rho(K_{2},u)\cdots\rho(K_{n},u)}{\varepsilon}~~~~~~~~~~~~~~~~~~~~~~~~~~~~~~~~~~~~~~~~~~~~~~~~~~~~~~~~~~~~~~~~$$
$$~~~=\rho(K_{2},u)\cdots\rho(K_{n},u)\lim_{\varepsilon\rightarrow 0^+}\frac{\rho_{\varepsilon}}{\varepsilon}\cdot\frac{\rho_{\varepsilon}-\rho(L_{1},u)}{\rho_{\varepsilon}}~~~~~~~~~~~~~~~~~~~~~~~~~~~~~~~~~~~~~$$
$$~~~~~~~~~~~~=\rho(K_{2},u)\cdots\rho(K_{n},u)\lim_{\varepsilon\rightarrow 0^+}\rho_{\varepsilon}\cdot\phi\left(\frac{\rho(K_{1},u)}{\rho_{\varepsilon}}\right)
\cdot\frac{1-\phi^{-1}\left(\phi(1)-\varepsilon{\phi}
\left(\frac{\displaystyle\rho(K_{1},u)}{\displaystyle\rho_{\varepsilon}}\right)\right)}{\phi(1)-\left(\phi(1)-\varepsilon
\phi\left(\frac{\displaystyle\rho(K_{1},u)}{\displaystyle\rho_{\varepsilon}}\right)\right)}.$$
Noting that $y\rightarrow 1^{+}$ as $\varepsilon\rightarrow
0^{+},$ we have
$$\lim_{\varepsilon\rightarrow 0^+}\frac{(\rho_{\varepsilon}-\rho(L_{1},u))\rho(K_{2},u)\cdots\rho(K_{n},u)}{\varepsilon}~~~~~~~~~~~~~~~~~~~~~~~~~~~~~~~~~~~~~~~~~~~~~~~~~~~~~~~~~~~~~~~~$$
$$=\rho(L_{1},u)\rho(K_{2},u)\cdots\rho(K_{n},u)\phi\left(\frac{\rho(K_{1},u)}{\rho(L_{1},u)}\right)\lim_{y\rightarrow 1^+}\frac{1-y}{\phi(1)-\phi(y)}~~~~~~~~~~~~~$$
$$=\frac{1}{\phi'_{r}(1)}\phi\left(\frac{\rho(K_{1},u)}{\rho(L_{1},u)}\right)\rho(L_{1},u)\rho(K_{2},u)\cdots\rho(K_{n},u),~~~~~~~~~~~~~~~~~~~~~~~~~~~\eqno(4.3)$$
where $$y=\phi^{-1}\left(\phi(1)-\varepsilon{\phi}
\left(\frac{\rho(K_{1},u)}{\rho_{\varepsilon}}\right)\right).$$
The equation (4.2) follows immediately from (2.3) with
(4.3).~~~~~~~~~~~~~~~~~~~~~~~~~~~~~~~~~~~~~~~~~~~~~~~~~~~~~~~~~~$\Box$

This theorem plays a central role in deriving the Orlicz multiple
dual mixed volume. Hence, we give the second proof.

{\bf Proof}~ Since
$$\frac{d\rho_{\varepsilon}}{d\varepsilon}=\frac{d}{d\varepsilon}\left(\frac{\displaystyle\rho(L_{1},u)}{\displaystyle{\phi}^{-1}
\left(\phi(1)-\varepsilon{\phi}
\left(\frac{\rho(K_{1},u)}{\rho_{\varepsilon}}\right)\right)}\right)~~~~~~~~~~~~~~~~~~~~~~~~~~~~~~~~~~~~~~~~~~~~~~~~~~~~~$$
$$=\frac{\rho(L_{1},u)\frac{\displaystyle d\phi^{-1}(y)}{\displaystyle dy}\left[{\phi}\left(\frac{\displaystyle
\rho(K_{1},u)}{\displaystyle
\rho_{\varepsilon}}\right)-\varepsilon\cdot\frac{\displaystyle
d\phi(z)}{\displaystyle
dz}\frac{\displaystyle\rho(K_{1},u)}{\displaystyle\rho_{\varepsilon}^{2}}\frac{\displaystyle
d\rho_{\varepsilon}}{\displaystyle d\varepsilon}\right]
}{\displaystyle\left(\displaystyle{\phi}^{-1}\left(\phi(1)-\varepsilon{\phi}
\left(\frac{\rho(K_{1},u)}{\rho_{\varepsilon}}\right)\right)\right)^{2}}.~~~~~~~~~~~~~~~~~~~~$$
where
$$y=\phi(1)-\varepsilon{\phi}
\left(\frac{\rho(K_{1},u)}{\rho_{\varepsilon}}\right),$$ and
$$z=\frac{\displaystyle\rho(K_{1},u)}{\displaystyle
\rho_{\varepsilon}}.$$ Hence
$$\frac{d\rho_{\varepsilon}}{d\varepsilon}=\frac{\rho(L_{1},u)\frac{\displaystyle d\phi^{-1}(y)}{\displaystyle dy}{\phi}\left(\frac{\displaystyle
\rho(K_{1},u)}{\displaystyle
\rho_{\varepsilon}}\right)}{\displaystyle\left(\displaystyle{\phi}^{-1}\left(\phi(1)-\varepsilon{\phi}
\left(\frac{\rho(K_{1},u)}{\rho_{\varepsilon}}\right)\right)\right)^{2}+\varepsilon\cdot\frac{\rho(K_{1},u)\rho(L_{1},u)}{\rho_{\varepsilon}^{2}}\frac{\displaystyle
d\phi^{-1}(y)}{\displaystyle dy}\frac{d\phi(z)}{dz}}.\eqno(4.4)$$
Since
$$\lim_{\varepsilon\rightarrow 0^+}\frac{(\rho_{\varepsilon}-\rho(L_{1},u))\rho(K_{2},u)\cdots\rho(K_{n},u)}{\varepsilon}~~~~~~~~~~~~~~~~~~~~~~~~~~~~~~~~~~~~~~~~~~~~~~~~~~~~~~~$$
$$=\rho(K_{2},u)\cdots\rho(K_{n},u)\lim_{\varepsilon\rightarrow 0^+}\frac{\rho_{\varepsilon}-\rho(L_{1},u)}{\varepsilon}~~~~~~~~~~~~~~~~~~~~~~~~~~~~~~~~~~~~~$$
$$=\rho(K_{2},u)\cdots\rho(K_{n},u)\lim_{\varepsilon\rightarrow 0^+}\frac{d\rho_{\varepsilon}}{d\varepsilon}.~~~~~~~~~~~~~~~~~~~~~~~~~~~~~~~~~~~~~~~~~~~~~~~~~\eqno(4.5)$$
On the other hand
$$\lim_{\varepsilon\rightarrow 0^+}\frac{\displaystyle
d\phi^{-1}(y)}{\displaystyle dy}=\lim_{\triangle y\rightarrow
0^+}\frac{\phi^{-1}(\phi(1)+\triangle y)-1}{\triangle
y}~~~~~~~~~~~$$
$$=\lim_{\omega\rightarrow 1^+}\frac{\omega-1}{\phi(\omega)-\phi(1)}~~~~~~$$
$$=\frac{1}{{\phi}^{'}_{r}(1)},~~~~~~~~~~~~~~~~~~~~~\eqno(4.6)$$
where $\omega=\phi^{-1}(\phi(1)+\triangle y).$

From (4.4), (4.5), (4.6) and Lemma 4.2, we obtain
$$\lim_{\varepsilon\rightarrow 0^+}\frac{(\rho_{\varepsilon}-\rho(L_{1},u))\rho(K_{2},u)\cdots\rho(K_{n},u)}{\varepsilon}=\frac{1}{{\phi}_{r}'(1)}\cdot{\phi}\left(\frac{\rho(K_{1},u)}{\rho(L_{1},u)}\right)
\rho(L_{1},u)\rho(K_{2},u)\cdots\rho(K_{n},u).\eqno(4.7)$$ From
(2.3) and (4.7), the equation (4.2) follows
easy.~~~~~~~~~~~~~~~~~~~~~~~~~~~~~~~~~~~~~~~~~~~~~~~~~~~~~~~~~~~~~~~~~~~~~~$\Box$

For any $L_{1},K_{1},\cdots,K_{n}\in {\cal S}^{n}$ and
${\phi}\in{\cal C}$, the integral on the right-hand side of (4.2)
denoting by $\widetilde{V}_{{\phi}}(L_{1},K_{1},\cdots,K_{n})$,
and hence this new Orlicz multiple dual mixed volume
$\widetilde{V}_{\phi}(L_{1},K_{1},\cdots ,K_{n})$ has been born.

{\bf Lemma 4.4}~ {\it If $L_{1},K_{1},\cdots,K_{n}\in {\cal
S}^{n}$ and $\phi\in{\cal C}$, then}
$$\widetilde{V}_{\phi}(L_{1},K_{1},\cdots ,K_{n})={\phi}'_{r}(1)\cdot\frac{d}{d\varepsilon}\bigg|_{\varepsilon=0^{+}}\widetilde{V}(L_{1}\widehat{+}_{\phi}\varepsilon\cdot
K_{1},K_{2},\cdots,K_{n}).\eqno(4.8)$$

{\bf Proof}~ This yields immediately from the Definition 4.1 and
the variational formula of volume (4.2)

{\bf Lemma 4.5}~ {\it Let $K,L\in {\cal S}^{n}$ and ${\phi}\in
{\cal C}$, then}

$$\lim_{\varepsilon\rightarrow 0^+}\frac{\widetilde{V}_{1}(K,K\widehat{+}_{\phi}\varepsilon\cdot
L)-V(K)}{\varepsilon}=\frac{1}{n}\lim_{\varepsilon\rightarrow
0^+}\frac{V(K\widehat{+}_{\phi}\varepsilon\cdot
L)-V(K)}{\varepsilon}.\eqno(4.9)$$

{\bf Proof}~  Suppose $\varepsilon>0$, $K,L\in {\cal S}^{n}$ and
$u\in S^{n-1}$, let
$$\bar{\rho}_{\varepsilon}=\rho(K\widehat{+}_{\phi}\varepsilon\cdot L,u).$$
From (1.3), (2.1), (2.2) and (4.4), we obtain
$$\lim_{\varepsilon\rightarrow 0^+}\frac{\widetilde{V}_{1}(K,K\widehat{+}_{\phi}\varepsilon\cdot
L)-V(K)}{\varepsilon}=\frac{1}{n}\int_{S^{n-1}}\lim_{\varepsilon\rightarrow
0^+}\frac{\rho(K\widehat{+}_{\phi}\varepsilon\cdot
L),u)\rho(K,u)^{n-1}-\rho(L,u)^{n}}{\varepsilon}dS(u)$$
$$~~~~~~~~~~~~~=\frac{1}{n}\int_{S^{n-1}}\rho(K,u)^{n-1}
\lim_{\varepsilon\rightarrow
0^+}\frac{d\bar{\rho}_{\varepsilon}}{d\varepsilon}dS(u)$$
$$~~~~~~~~~~~~~~~~~~~~~~=\frac{1}{n\phi'_{r}(1)}\int_{S^{n-1}}\phi\left(\frac{\rho(L,u)}{\rho(K,u)}\right)
\rho(K,u)^{n}dS(u)$$
$$~~~~~~~~~~~~~~~~~~~~~~~~~~~~~~~~~~~~~~~~~~~~~~~~=\frac{1}{n}\lim_{\varepsilon\rightarrow
0^+}\frac{V(K\widehat{+}_{\phi}\varepsilon\cdot
L)-V(K)}{\varepsilon}.~~~~~~~~~~~~~~~~~~~~~~~~~~~~~~~~~~~~~~\Box$$

{\bf Lemma 4.6}~ {\it Let $L_{1},K_{1},\ldots,K_{n}\in {\cal
S}^{n}$ and ${\phi}\in {\cal C}$, then}
$$\widetilde{V}_{\phi}(L_{1},K_{1},\cdots
,K_{n})=\widetilde{V}_{\phi}(K,L),\eqno(4.10)$$ if
$K_{2}=\cdots=K_{n}=K$, $L_{1}=K$ and $K_{1}=L.$

{\bf Proof}~ On the one hand, putting $K_{2}=\cdots=K_{n}=K$,
$L_{1}=K$ and $K_{1}=L$ in (4.8), and noting Lemma 4.5 and (1.3),
it follows that
$$\widetilde{V}_{\phi}(L_{1},K_{1},\cdots ,K_{n})=\phi'_{r}(1)\frac{d}{d\varepsilon}\bigg|_{\varepsilon=0^{+}}\widetilde{V}(L_{1}\widehat{+}_{\phi}\varepsilon\cdot
K_{1},K_{2},\cdots,K_{n})~~~~~~~~~~~~~~~~~~~~~~~~~~$$
$$=\phi'_{r}(1)\lim_{\varepsilon\rightarrow 0^+}\frac{\widetilde{V}_{1}(K,K\widehat{+}_{\phi}\varepsilon\cdot
L)-V(K)}{\varepsilon}~~~~~$$
$$=\frac{\phi'_{r}(1)
}{n}\lim_{\varepsilon\rightarrow
0^+}\frac{V(K\widehat{+}_{\phi}\varepsilon\cdot
L)-V(K)}{\varepsilon}~~~~~~~~~$$
$$=\widetilde{V}_{\phi}(K,L).~~~~~~~~~~~~~~~~~~~~~~~~~~~~~~~~~~~~~~~~~\eqno(4.11)$$
On the other hand, let $K_{2}=\cdots=K_{n}=K$, $L_{1}=K$ and
$K_{1}=L$, from Definition 4.1 and (1.3), then
$$\widetilde{V}_{\phi}(L_{1},K_{1},\cdots
,K_{n})=\frac{1}{n}\int_{S^{n-1}}{\phi}
\left(\frac{\rho(K_{1},u)}{\rho(L_{1},u)}\right)\rho(L_{1},u)\rho(K_{2},u)\cdots\rho(K_{n},u)
dS(u)~~~~$$
$$=\frac{1}{n}\int_{S^{n-1}}\phi\left(\frac{\rho(L,u)}{\rho(K,u)}\right)
\rho(K,u)^{n}dS(u)~~~~~$$
$$=\widetilde{V}_{\phi}(K,L).~~~~~~~~~~~~~~~~~~~~~~~~~~~~~~~~~~~~~~~~\eqno(4.12)$$
Combining (4.11) and (4.12), this shows that
$$\widetilde{V}_{\phi}(L_{1},K_{1},\cdots
,K_{n})=\widetilde{V}_{\phi}(K,L),$$ if $K_{2}=\cdots=K_{n}=K$,
$L_{1}=K$ and $K_{1}=L$.
~~~~~~~~~~~~~~~~~~~~~~~~~~~~~~~~~~~~~~~~~~~~~~~~~~~~~~~~~~~~~~~~~~~~~~~~~~~~$\Box$

{\bf Lemma 4.7}~ [60] {\it If $K_{i},L_{i}\in {\cal S}^{n}$ and
$K_{i}\rightarrow K$, $L_{i}\rightarrow L$ as $
i\rightarrow\infty$, then
$$a\cdot K_{i}\widehat{+}_{\phi}b\cdot L_{i}\rightarrow a\cdot K\widehat{+}_{\phi}b\cdot L,~~ {as}~~ i\rightarrow\infty,\eqno(4.13)$$
for all $a$ and $b$.}

{\bf Lemma 4.8}~ {\it If $L_{1},K_{1},\ldots,K_{n},K,L\in {\cal
S}^{n}$, $\lambda_{1},\cdots,\lambda_{n}\geq 0$ and ${\phi}\in
{\cal C}$, then}

(1) $\widetilde{V}_{\phi}(L_{1},K_{1},\cdots ,K_{n})>0$

(2)
$\widetilde{V}_{\phi}(K_{1},K_{1},K_{2},\cdots,K_{n})=\phi(1)\widetilde{V}(K_{1},\cdots,K_{n})$

(3) $\widetilde{V}_{\phi}(K,K,\cdots,K)=\phi(1)V(K)$

(4)
$$\widetilde{V}_{\phi}(\lambda_{1}L_{1},\lambda_{1}K_{1},\lambda_{2}K_{2}\cdots,\lambda_{n}K_{n})=\lambda_{1}\cdots\lambda_{n}\widetilde{V}_{\phi}(L_{1},K_{1},\cdots,K_{n}).$$

(5) $$\widetilde{V}_{\phi}(L_{1}, K_{1},\lambda_{1}
K\widetilde{+}\lambda_{2}
L,K_{3}\cdots,K_{n})~~~~~~~~~~~~~~~~~~~~~~~~~~~~~~~~~~~~~~~~~~~~~~~~~~~~~~~~~~~~~~~~~~~~$$
$$=\lambda_{1}
V_{\phi}(L_{1},K_{1},K,K_{3},\cdots,K_{n})+\lambda_{2}
V_{\phi}(L_{1},K_{1},L,K_{3},\cdots,K_{n}).$$

This shows the Orlicz multiple mixed volume
$V_{\phi}(K_{1},\cdots,K_{n},L_{n})$ is linear in its back $(n-1)$
variables.

(6) $\widetilde{V}_{\phi}(L_{1},K_{1},\cdots ,K_{n})$ is
continuous.

{\bf Proof}~ From Definition 4.1, it immediately gives (1), (2),
(3) and (4).

From Definition 4.1, combining the following fact
$$\rho(\lambda_{1}K\widetilde{+}\lambda_{2} L,\cdot)=\lambda_{1}\rho(K,\cdot)+\lambda_{2}\rho(L,\cdot),$$
it yields (5) directly.

Suppose $L_{i1}\rightarrow L_{1}$, $K_{ij}\rightarrow K_{j}$ as
$i\rightarrow\infty$ where $j=1,\ldots,n$, combining Definition
4.1 and Lemma 4.6 with the following facts
$$\widetilde{V}(L_{i1}\widehat{+}_{\phi}\varepsilon\cdot
K_{i1},K_{i2},\cdots,K_{in})\rightarrow
\widetilde{V}(L_{1}\widehat{+}_{\phi}\varepsilon\cdot
K_{1},K_{2},\cdots,K_{n})$$ and
$$\widetilde{V}(L_{i1},K_{i2},\cdots,K_{in})\rightarrow\widetilde{V}(L_{1},K_{2},\cdots,K_{n})$$
as $i\rightarrow\infty,$ it yield (6) directly.
~~~~~~~~~~~~~~~~~~~~~~~~~~~~~~~~~~~~~~~~~~~~~~~~~~~~~~~~~~~~~~~~~~~~~~~~~~~~~~~~~~~~~~~~~~~~~~~~~~~~~~~~~~~$\Box$

{\bf Lemma 4.9}~ [60] {\it Suppose $K,L\in{\cal S}^{n}$ and
$\varepsilon>0$. If $\phi\in {\cal C}$, then for} $A\in GL(n)$
$$A(K\widetilde{+}_{\phi}\varepsilon\cdot L)=AK\widetilde{+}_{\phi}\varepsilon\cdot AL.\eqno(4.14)$$

We easy find that Orlicz multiple dual mixed volume
$\widetilde{V}_{\phi}(L_{1},K_{1},K_{2},\cdots ,K_{n})$ is
invariant under simultaneous unimodular centro-affine
transformation.

{\bf Lemma 4.10}~ {\it If $L_{1},K_{1},\ldots,K_{n}\in {\cal
S}^{n}$ and ${\phi}\in {\cal C}$, then for $A\in SL(n)$,}
$$\widetilde{V}_{\phi}(AL_{1},AK_{1},\cdots,AK_{n})=\widetilde{V}_{\phi}(L_{1},K_{1},\cdots,K_{n}).\eqno(4.15)$$

{\bf Proof}~ From (4.8) and Lemma 4.8, we have, for $A\in SL(n)$,
$$\widetilde{V}_{\phi}(AL_{1},AK_{1},\cdots,AK_{n})~~~~~~~~~~~~~~~~~~~~~~~~~~~~~~~~~~~~~~~~~~~~~~~~~~~~~~~~~~~~~~~~~~~~~~~~~~~$$
$$={\phi}'_{r}(1)\lim_{\varepsilon\rightarrow 0^+}
\frac{\widetilde{V}(AL_{1}\widehat{+}_{\phi}\varepsilon\cdot
AK_{1},AK_{2},\cdots,AK_{n})-\widetilde{V}(AL_{1},AK_{2},\cdots,AK_{n})}{\varepsilon}$$
$$={\phi}'_{r}(1)\lim_{\varepsilon\rightarrow 0^+}
\frac{\widetilde{V}(A(L_{1}\widehat{+}_{\phi}\varepsilon\cdot
K_{1}),AK_{2},\cdots,AK_{n})-\widetilde{V}(L_{1},K_{2},\cdots,K_{n})}{\varepsilon}~~~~~~~$$
$$={\phi}'_{r}(1)\lim_{\varepsilon\rightarrow 0^+}
\frac{\widetilde{V}(L_{1}\widehat{+}_{\phi}\varepsilon\cdot
K_{1},K_{2},\cdots,K_{n})-\widetilde{V}(L_{1},K_{2},\cdots,K_{n})}{\varepsilon}~~~~~~~~~~~~~~~~$$
$$=\widetilde{V}_{\phi}(L_{1},K_{1},\cdots,K_{n}).~~~~~~~~~~~~~~~~~~~~~~~~~~~~~~~~~~~~~~~~~~~~~~~~~~~
~~~~~~~~~~~~~~~\Box$$

For the convenience of writing, when $K_{1}=\cdots=K_{i}=K$,
$K_{i+1}=\cdots=K_{n}=L$, $L_{n}=M$, the Orlicz multiple dual
mixed volume $\widetilde{V}_{\phi}(K,\cdots,K,L,\cdots,L,M)$, with
$i$ copies of $K$, $n-i$ copies of $L$ and $1$ copy of $M$, will
be denoted by $\widetilde{V}_{\phi}(K~[i],L~[n-i],M)$.

{\bf Lemma 4.11}~ {\it If $K,L\in {\cal S}^{n}$ and ${\phi}\in
{\cal C}$, and $0\leq i<n$ then}
$$\widetilde{V}_{\phi}(K,L,K~[n-i-1],B~[i])=\frac{1}{n}\int_{S^{n-1}}\phi
\left(\frac{\rho(L,u)}{\rho(K,u)}\right)\rho(K,u)^{n-i}dS(u).\eqno(4.16)$$

{\bf Proof}~ On the one hand, putting $L_{1}=K$,
$K_{1}=L$,$K_{2}=\cdots=K_{n-i}=K$ and $K_{n-i+1}=\cdots=K_{n}=B$
in (4.8), from (2.4), (2.5), (4.4) and (4.6), we obtain for
$\varphi_{1}, \varphi_{2}\in \Phi$
$$~~~~~~~~~~~~~~~~~~~~~\widetilde{V}_{\phi}(K,L,K~[n-i-1],B~[i])=\phi^{'}_{r}(1)\lim_{\varepsilon\rightarrow 0^+}\frac{\widetilde{W}_{i}(K,K\widehat{+}_{\phi
}\varepsilon\cdot
L)-\widetilde{W}_{i}(K)}{\varepsilon}~~~~~~~~~~~~~~~~~~~~~~~~~~~~~~$$
$$~~~~~~~~~~~~~~~~~~~~~~~~~~~~~=\frac{1}{n}\phi^{'}_{r}(1)\int_{S^{n-1}}\lim_{\varepsilon\rightarrow 0^+}\frac{\rho(K\widehat{+}_{\phi}\varepsilon\cdot L)\rho(K,u)^{n-i-1}-\rho(K,u)^{n-i}}{\varepsilon}dS(
u)$$
$$=\frac{1}{n}\phi^{'}_{l}(1)\int_{S^{n-1}}\rho(K,u)^{n-i-1}\lim_{\varepsilon\rightarrow 0^+}\frac{d\overline{\rho}_{\varepsilon}}{\varepsilon}dS(u)$$
$$=\frac{1}{n}\int_{S^{n-1}}\phi
\left(\frac{\rho(L,u)}{\rho(K,u)}\right)\rho(K,u)^{n-i}dS(u).~~~~\eqno(4.17)$$
On the one hand, putting $L_{1}=K$,
$K_{1}=L$,$K_{2}=\cdots=K_{n-i}=K$ and $K_{n-i+1}=\cdots=K_{n}=B$
in Definition 4.1, we have
$$\widetilde{V}_{\phi}(K,L,K~[n-i-1],B~[i])=\frac{1}{n}\int_{S^{n-1}}\phi
\left(\frac{\rho(L,u)}{\rho(K,u)}\right)\rho(K,u)^{n-i}dS(u).\eqno(4.18)$$
Combining (4.17) and (4.18), (4.16) yields
easy.~~~~~~~~~~~~~~~~~~~~~~~~~~~~~~~~~~~~~~~~~~~~~~~~~~~~~~~~~~~~~~~~~~~~~~~~~~~~~~~~~~~~~~~~~~~~~~~~~~~~~~$\Box$

Here, we denote the Orlicz multiple dual mixed volume
$\widetilde{V}_{\phi}(K,L,K~[n-i-1],B~[i])$ by
$\widetilde{W}_{\phi,i}(K,L)$, and call
$\widetilde{W}_{\phi,i}(K,L)$ as Orlicz dual quermassintegral of
star bodies $K$ and $L$. When $i=0$, Orlicz dual quermassintegral
$\widetilde{W}_{\phi,i}(K,L)$ becomes Orlicz dual mixed volume
$\widetilde{V}_{\phi}(K,L)$.

{\bf Remark 4.12}~ When $\phi(t)=t^{-p}$, $p=1$ and $L_1=K_1$,
from the integral representation for Orlicz multiple dual mixed
volume $\widetilde{V}_{\phi}(L_{1},K_{1},\cdots ,K_{n})$ (see
Definition 4.1), it follows easy that
$\widetilde{V}_{\phi}(L_{1},K_{1},\cdots
,K_{n})=\widetilde{V}(K_{1},\cdots,K_{n}).$ On the other hand,
when $\phi(t)=t^{-p}$, $p=1$ and $L_1=K_1$, from (4.8) and noting
that $\phi'_{r}(1)=-1$, hence
$$\widetilde{V}(K_{1},\cdots,K_{n})=\lim_{\varepsilon\rightarrow 0^+}
\frac{\widetilde{V}(K_{1},\cdots
,K_{n})-\widetilde{V}(K_{1}\widehat{+}\varepsilon\cdot
K_{1},K_{2},\cdots,K_{n})}{\varepsilon}.\eqno(4.19)$$ This is very
interesting for the usually dual mixed volume of this form.

{\bf Remark 4.13}~ When $\phi(t)=t^{-p}$, $p\geq 1$, write the
Orlicz multiple dual mixed volume
$\widetilde{V}_{\phi}(L_{1},K_{1},\cdots ,K_{n})$ as
$\widetilde{V}_{-p}(L_{1},K_{1},\cdots ,K_{n})$ and call as the
$L_p$-multiple dual mixed volume, from Definition 4.1, it easy
yields
$$\widetilde{V}_{-p}(L_{1},K_{1},\cdots,K_{n})=\frac{1}{n}\int_{S^{n-1}}\rho(K_{1},u)^{-p}\rho(L_{1},u)^{1+p}\rho(K_{2},u)\cdots\rho(K_{n},u)
dS(u).$$ When $\varphi(t)=t^{-p}$ and $p\geq 1$, from (4.8), we
get the following expression of $L_{p}$-multiple dual mixed
volume.
$$\frac{1}{-p}\widetilde{V}_{-p}(L_{1},K_{1},\cdots ,K_{n})=\lim_{\varepsilon\rightarrow 0^+}
\frac{\widetilde{V}(L_{1}\widehat{+}_{p}\varepsilon\cdot
K_{1},K_{2},\cdots,K_{n})-\widetilde{V}(L_{1},K_{2},\cdots
,K_{n})}{\varepsilon}.$$

{\bf Lemma 4.14}~ (Jensen's inequality) {\it Let $\mu$ be a
probability measure on a space $X$ and $g: X\rightarrow I\subset
{\Bbb R}$ is a $\mu$-integrable function, where $I$ is a possibly
infinite interval. If ${\psi}: I\rightarrow {\Bbb R}$ is a convex
function, then
$$\int_{X}{\psi}(g(x))d\mu(x)\geq{\psi}\left(\int_{X}g(x)d\mu(x)\right).\eqno(4.20)$$ If ${\psi}$ is strictly convex, equality holds if and only if
$g(x)$ is constant for $\mu$-almost all $x\in X$} (see [24, p.165]).

\vskip 10pt \noindent{\large \bf 5 ~ Dual Orlicz-Aleksandrov-Fenchel inequality}\vskip 10pt

{\bf Theorem 5.1}~ {\it If $L_{1},K_{1},\cdots,K_{n}\in {\cal
S}^{n}$ and ${\phi}\in {\cal C}$, then
$$\widetilde{V}_{\phi}(L_{1},K_{1},\cdots,
K_{n})\geq\widetilde{V}(L_{1},K_{2},\cdots, K_{n})
\phi\left(\frac{\widetilde{V}(K_{1},\cdots,K_{n})}{\widetilde{V}(L_{1},K_{2},\cdots,K_{n})}\right).\eqno(5.1)$$
If $\phi$ is strictly convex, equality holds if and only if
$K_{1}$ and $L_{1}$ are dilates.}

{\bf Proof}~ For $K_{1},\cdots,K_{n}\in {\cal S}^{n}$ and any
$u\in S^{n-1}$, since
$$\frac{1}{n\widetilde{V}(K_{1},\cdots,K_{n})}\int_{S^{n-1}}\rho(K_{1},u)\cdots\rho(K_{n},u)dS(u)=1,$$
so
$\frac{\displaystyle\rho(K_{1},u)\cdots\rho(K_{n},u)}{\displaystyle
n\widetilde{V}(K_{1},\cdots,K_{n})}S(u)$ is a probability measure
on $S^{n-1}$. From Definition 4.1, Jensen's inequality (4.20) and
(2.3), it follows that
$$\frac{\widetilde{V}_{\phi}(L_{1},K_{1},\cdots,K_{n})}{\widetilde{V}(L_{1},K_{2},\cdots,K_{n})}=\frac{1}{n\widetilde{V}(L_{1},K_{2},\cdots,K_{n})}~~~~~~~~~~~~~~~~~~~~~~~~~~~~~~~~~~~~~~~~~~~~~~~~~~~~~~~~~~~~~~~~~~~~~~~~~~$$
$$~~~~~~~\times\int_{S^{n-1}}{\phi}
\left(\frac{\rho(K_{1},u)}{\rho(L_{1},u)}\right)\rho(L_{1},u)\rho(K_{2},u)\cdots\rho(K_{n},u)dS(u)$$
$$\geq\phi\left(\frac{1}{n\widetilde{V}(L_{1},K_{2},\cdots,K_{n})}\int_{S^{n-1}}\rho(K_{1},u)
\cdots\rho(K_{n},u)dS(u)\right)$$
$$=\phi\left(\frac{\widetilde{V}(K_{1},\cdots,K_{n})}{\widetilde{V}(L_{1},K_{2},\cdots,K_{n})}\right).~~~~~~~~~~~~~~~~~~~~~~~~~~~~~~~~~~~~~~~~~~~~\eqno(5.2)$$
Next, we discuss the equality condition of (5.2). Suppose the
equality hold in (5.2), form the equality condition of Jensen's
inequality, it follows that if $\phi$ is strictly convex the
equality in (5.2) holds if and only if $K_{1}$ and $L_{1}$ are
dilates.
~~~~~~~~~~~~~~~~~~~~~~~~~~~~~~~~~~~~~~~~~~~~~~~~~~~~~~~~~~~~~~~~~~~~~~~~~~~~~~~~~~~~~~~~~~~~~~~~~~~~~~~~$\Box$

{\bf Theorem 5.2}~ (Dual Orlicz-Aleksandrov-Fenchel inequality)
{\it If $L_{1},K_{1},\cdots,K_{n}\in {\cal S}^{n}$, $1\leq r\leq
n$ and ${\phi}\in {\cal C}$, then
$$\widetilde{V}_{\phi}(L_{1},K_{1},K_{2},\cdots,
K_{n})\geq\widetilde{V}(L_{1},K_{2},\cdots, K_{n})
\phi\left(\frac{\prod_{i=1}^{r}\widetilde{V}(K_{i}\ldots,K_{i},K_{r+1},\ldots,K_{n})^{\frac{1}{r}}}{\widetilde{V}(L_{1},
K_{2}\ldots,K_{n})}\right).\eqno(5.3)$$ If $\phi$ is strictly
convex, equality holds if and only if $L_{1}, K_{1},\ldots,K_{r}$
are all dilations of each other.}

{\bf Proof}~ This follows immediately from Theorem 5.1 with the
dual Aleksandrov-Fenchel inequality. ~~~$\Box$

Obviously, putting $\phi(t)=t^{-p}$, $p=1$ and $L_1=K_1$ in (5.3),
(5.3) becomes the Lutwak's dual Aleksandrov-Fenchel inequality
(1.11) stated in the introduction.

{\bf Corollary 5.3}~ {\it If $L_{1},K_{1},\cdots,K_{n}\in {\cal
S}^{n}$ and ${\phi}\in {\cal C}$, then
$$\widetilde{V}_{\phi}(L_{1},K_{1},K_{2},\cdots,
K_{n})\geq\widetilde{V}(L_{1},K_{2},\cdots, K_{n})
\phi\left(\left(\frac{V(K_{1})\cdots
V(K_{n})}{\widetilde{V}(L_{1},K_{2},\ldots,K_{n})^{n}}\right)^{\frac{1}{n}}\right).\eqno(5.4)$$
If $\phi$ is strictly convex, equality holds if and only if
$L_{1}, K_{1},\ldots,K_{n}$ are all dilations of each other.}

{\bf Proof}~ This follows immediately from Theorem 5.2 with $r=n$.
~~~~~~~~~~~~~~~~~~~~~~~~~~~~~~~~~~~~~~~~~~~~~~~~~~~~~~~~~~~~$\Box$

{\bf Corollary 5.4}~ {\it If $K,L\in {\cal S}^{n}$, $0\leq i<n$
and $\phi\in {\cal C}$, then
$$\widetilde{W}_{\phi,i}(K,L)\geq
\widetilde{W}_{i}(K)\phi\left(\left(\frac{\widetilde{W}_{i}(L)}{\widetilde{W}_{i}(K)}\right)^{1/(n-i)}\right)
.\eqno(5.5)$$ If $\phi$ is strictly convex, equality holds if and
only if $K$ and $L$ are dilates.}

{\bf Proof}~ This follows immediately from Theorem 5.2 with
$r=n-i$, $L_{1}=K$, $K_{1}=L$, $K_{2}=\cdots=K_{n-i}=K$ and
$K_{n-i+1}=\cdots=K_{n}=B$.
~~~~~~~~~~~~~~~~~~~~~~~~~~~~~~~~~~~~~~~~~~~~~~~~~~~~~~~~~~~~~~~~~~~~~~~~~~~~~~~~~~~~~~~~~~~~~~~~~~~~~~~~~~~$\Box$

The following inequality follows immediately from (5.5) with
$\phi(t)=t^{-p}$ and $p\geq 1$. If $K,L\in {\cal S}^{n}$, $0\leq
i<n$ and $p\geq 1$, then
$$\widetilde{W}_{-p,i}(K,L)^{n-i}\geq\widetilde{W}_{i}(K)^{n-i+p}\widetilde{W}_{i}(L)^{-p},\eqno(5.6)$$
with equality if and only if $K$ and $L$ are dilates. Taking $i=0$
in (5.6), this yields Lutwak's $L_{p}$-dual Minkowski inequality
is following: If $K,L\in {\cal S}^{n}$ and $p\geq 1$, then
$$\widetilde{V}_{-p}(K,L)^{n}\geq V(K)^{n+p} V(L)^{-p},\eqno(5.7)$$
with equality if and only if $K$ and $L$ are dilates.

{\bf Theorem 5.5} (Orlicz dual isoperimetric inequality) {\it If
$K\in {\cal S}^{n}$ and $\phi\in{\cal C}$, and $0\leq i<n$ then
$$\frac{\widetilde{V}_{\phi}(K,B,K~[n-i-1],B~[i])}{\widetilde{W}_{i}(K)}\geq
\phi\left(\left(\frac{V(B)}{\widetilde{W}_{i}(K)}\right)^{1/(n-i)}\right).\eqno(5.8)$$
If $\phi$ is strictly convex, equality holds if and only if $K$ is
a ball.}

{\bf Proof}~ This follows immediately from (5.3) with with
$r=n-i$, $L_{1}=K$, $K_{1}=B$, $K_{2}\cdots=K_{n-i}=K,$ and $
K_{n-i+1}=\cdots=K_{n}=B$.
~~~~~~~~~~~~~~~~~~~~~~~~~~~~~~~~~~~~~~~~~~~~~~~~~~~~~~~~~~~~~~~~~~~~~~~~~~~~~~~~~~~~~~~~~~~~~~~~~~~~~~~~~~~~~$\Box$

When $\phi(t)=t^{-p},$ $p\geq 1$, the Orlicz isoperimetric
inequality (5.8) becomes the following $L_p$-dual isoperimetric
inequality. If $K$ is a star body, $p\geq 1$ and $0\leq i<n$, then
$$\left(\frac{n\widetilde{V}_{-p}(K,B)}{\omega_{n}}\right)^{n-i}\geq\left(\frac{\widetilde{W}_{i}(K)}{\kappa_{n}}\right)^{n-i+p},\eqno(5.9)$$
with equality if and only if $K$ is ball, and where $\kappa_{n}$
denotes volume of the unit ball $B$, and its surface area by
$\omega_{n}.$

Putting $p=1$ and $i=0$ in (5.9), (5.9) becomes the following dual
isoperimetric inequality. If $K$ is a star body, then
$$\left(\frac{n\widetilde{V}_{-1}(K,B)}{\omega_{n}}\right)^{n}\geq\left(\frac{V(K)}{\kappa_{n}}\right)^{n+1},$$
with equality if and only if $K$ is ball.

{\bf Theorem 5.6}~ {\it If $L_{1},K_{1},\cdots,K_{n}\in {\cal
M}\subset{\cal S}^{n}$, and $\phi\in{\cal C}$ be strictly convex,
and if either
$$\widetilde{V}_{\phi}(Q,K_{1},\ldots,K_{n})=\widetilde{V}_{\phi}(Q,L_{1},K_{2},\ldots,K_{n}),~~ {for~ all}~~ Q\in{\cal M,}\eqno(5.10)$$
or
$$\frac{\widetilde{V}_{\phi}(K_{1},Q,K_{2},\ldots,K_{n})}{\widetilde{V}(K_{1},\ldots,K_{n})}=\frac{\widetilde{V}_{\phi}(L_{1},Q,K_{2},\ldots,K_{n})}{\widetilde{V}(L_{1},K_{2},\ldots,K_{n})},~~ {for~ all}~~ Q\in{\cal M,}\eqno(5.11)$$
then} $K_{1}=L_{1}.$

{\bf Proof}~ Suppose (5.10) hold. Taking $K_{1}$ for $Q$, then
from Definition 4.1 and Theorem 5.1, we obtain
$$\phi(1)\widetilde{V}(K_{1},\ldots,K_{n})=\widetilde{V}_{\phi}(K_{1},L_{1},K_{2},\ldots,K_{n})
\geq
\widetilde{V}(K_{1},\ldots,K_{n})\phi\left(\frac{\widetilde{V}(L_{1},K_{2},\ldots,K_{n})}{\widetilde{V}(K_{1},\ldots,K_{n})}\right),$$
with equality if and only if $K_{1}$ and $L_{1}$ are dilates.
Hence
$$\phi(1)\geq\phi\left(\frac{\widetilde{V}(L_{1},K_{2},\ldots,K_{n})}{\widetilde{V}(K_{1},\ldots,K_{n})}\right),$$
with equality if and only if $K_{1}$ and $L_{1}$ are dilates.
Since $\varphi$ is decreasing function on $(0,\infty),$ this
follows that
$$\widetilde{V}(K_{1},\ldots,K_{n})\leq \widetilde{V}(L_{1},K_{2},\ldots,K_{n}),$$
with equality if and only if $K_{1}$ and $L_{1}$ are dilates. On
the other hand, if taking $L_{1}$ for $Q$, we similar get
$\widetilde{V}(K_{1},\ldots,K_{n})\geq
\widetilde{V}(L_{1},K_{2},\ldots,K_{n}),$ with equality if and
only if $K_{1}$ and $L_{1}$ are dilates. Hence
$\widetilde{V}(K_{1},\ldots,K_{n})=\widetilde{V}(L_{1},K_{2},\ldots,K_{n}),$
and $K_{1}$ and $L_{1}$ are dilates, it follows that $K_{1}$ and
$L_{1}$ must be equal.

Suppose (5.11) hold. Taking $K_{1}$ for $Q$, then from Definition
4.1 and Theorem 5.1, we obtain
$$\phi(1)=\frac{\widetilde{V}_{\phi}(L_{1},K_{1},\ldots,K_{n})}{\widetilde{V}(L_{1},K_{2},\ldots,K_{n})}
\geq
\phi\left(\frac{\widetilde{V}(K_{1},\ldots,K_{n})}{\widetilde{V}(L_{1},K_{2}\ldots,K_{n})}\right),$$
with equality if and only if $K_{1}$ and $L_{1}$ are dilates.
Since $\varphi$ is increasing function on $(0,\infty),$ this
follows that
$$\widetilde{V}(L_{1},K_{2},\ldots,K_{n})\leq \widetilde{V}(K_{1},\ldots,K_{n}),$$
with equality if and only if $K_{1}$ and $L_{1}$ are dilates. On
the other hand, if taking $L_{1}$ for $Q$, we similar get
$\widetilde{V}(L_{1},K_{2},\ldots,K_{n})\geq
\widetilde{V}(K_{1},\ldots,K_{n}),$ with equality if and only if
$K_{1}$ and $L_{1}$ are dilates. Hence
$\widetilde{V}(L_{1},K_{2},\ldots,K_{n})=\widetilde{
V}(K_{1},\ldots,K_{n})$, and $K_{1}$ and $L_{1}$ are dilates, it
follows that $K_{1}$ and $L_{1}$ must be equal.
~~~~~~~~~~~~~~~~~~~~~~~~~~~~~~~~~~~~~~~~~~~~~~~~~~~~~~~~~~~~~~~~~~~~~~~~~~~~~~~~~~~~~~~~~~~~~~~~~~~~~~~~~~~~~~~~~~~~~~~~~~~~~~~~~~~~~~~~~~~~~~~~$\Box$

{\bf Corollary 5.7} {\it Let $K,L\in {\cal M}\subset{\cal S}^{n}$,
$0\leq i<n$, and $\phi\in {\cal C}$ be strictly convex, and if
either
$$\widetilde{W}_{\phi,i}(Q,K)=\widetilde{W}_{\phi,i}(Q,L),~ {\it for~ all}~ Q\in{\cal M},$$
or
$$\frac{\widetilde{W}_{\phi,i}(K,Q)}{\widetilde{W}_{i}(K)}=\frac{\widetilde{W}_{\phi,i}(L,Q)}{\widetilde{W}_{i}(L)},
~ {\it for~ all}~ Q\in{\cal M},$$ then} $K=L.$

{\bf Proof}~ This yields immediately from Theorem 5.6 and Lemma
4.11.~~~~~~~~~~~~~~~~~~~~~~~~~~~~~~~~~~~~~~~~~~~~~~~~~~~~~~~~~~~~~~~~~~$\Box$

{\bf Remark 5.8} When $\phi(t)=t^{-p}$ and $p=1$, the Orlicz dual
Aleksandrov-Fenchel inequality (5.3) becomes the following
inequality.  If $L_{1},K_{1},\cdots,K_{n}\in {\cal S}^{n}$ and
$1\leq r\leq n$, then
$$\widetilde{V}_{-1}(L_{1},K_{1},K_{2},\cdots,
K_{n})\geq\frac{\widetilde{V}(L_{1},K_{2},\cdots,
K_{n})^{2}}{\displaystyle\prod_{i=1}^{r}\widetilde{V}(K_{i}\ldots,K_{i},K_{r+1},\ldots,K_{n})^{\frac{1}{r}}},\eqno(5.12)$$
with equality if and only if $L_{1},K_{1},\ldots,K_{r}$ are all
dilations of each other.

Putting $L_{1}=K_{1}$ in (5.12) and noting that
$\widetilde{V}_{-1}(K_{1},K_{1},K_{2},\cdots,
K_{n})=\widetilde{V}(K_{1},K_{2},\cdots, K_{n})$, (5.12) becomes
the dual Aleksandrov-Fenchel inequality (1.11). Putting $r=n$ in
(5.12), (5.12) becomes the following inequality.
$$\widetilde{V}_{-1}(L_{1},K_{1},\cdots,
K_{n})^{n}\geq\widetilde{V}(L_{1},K_{2}\cdots,
K_{n})^{2n}(V(K_{1})\cdots V(K_{n}))^{-1},\eqno(5.13)$$ with
equality if and only if $L_{1},K_{1},\ldots,K_{n}$ are all
dilations of each other. Putting $L_{1}=K$, $K_{1}=L$ and
$K_{2}=\cdots=K_{n}=K$ in (5.13), (5.13) becomes the well-known
Minkowski inequality. If $K,L\in {\cal S}^{n}$, then
$$\widetilde{V}_{-1}(K,L)^{n}\geq V(K)^{n+1}V(L)^{-1},\eqno(5.14)$$
with equality if and only if $K$ and $L$ are dilates. Obviously,
inequality (5.12) in special case yields also the following
result. If $K_{1},\ldots,K_{n}\in {\cal S}^{n}$ and $0\leq i<n$,
then
$$\widetilde{W}_{-1,i}(K,L)^{n-i}\geq \widetilde{W}_{i}(K)^{n-i+1}W_{i}(L)^{-1},\eqno(5.15)$$
with equality if and only if $K$ and $L$ are dilates. When $i=0$,
(5.15) becomes (5.14). On the other hand, putting $L_{1}=K_{1}$ in
(5.13), (5.13) becomes the well-known inequality. If
$K_{1},\ldots,K_{n}\in {\cal S}^{n}$, then
$$\widetilde{V}(K_{1},\cdots,
K_{n})^{n}\leq V(K_{1})\cdots V(K_{n}),$$ with equality if and
only if $K_{1},\ldots,K_{n}$ are all dilations of each other.

\vskip 10pt \noindent{\large \bf 6 ~ Dual Orlicz-Brunn-Minkowski
inequality for the Orlicz harmonic addition}\vskip 10pt

{\bf Lemma 6.1}~ {\it If $L_{1},K_{1},\cdots,K_{n}\in {\cal
S}^{n}$ and ${\phi}\in {\cal C}$, then}
$$\phi(1)\widetilde{V}(K_{1}\widehat{+}_{{\phi}}L_{1},K_{2},\cdots,K_{n})=\widetilde{V}_{\phi}(K_{1}\widehat{+}_{\phi}L_{1},K_{1},\cdots,K_{n})+\widetilde{V}_{\phi}
(K_{1}\widehat{+}_{\phi}L_{1},L_{1},K_{2},\cdots,
K_{n}).\eqno(6.1)$$

{\bf Proof}~  Suppose $\varepsilon>0$, $K_{1},L_{1}\in {\cal
S}^{n}$ and $u\in S^{n-1}$, let
$$Q=K_{1}\widehat{+}_{\phi}\varepsilon\cdot L_{1}.$$
From Definition 4.1, (2.3)  and (3.4), we have
$$\phi(1)\widetilde{V}(Q,K_{2},\cdots,K_{n})~~~~~~~~~~~~~~~~~~~~~~~~~~~~~~~~~~~~~~~~~~~~~~~~~~~~~~~~~~~~~~~~~~~~~~~~~~~~$$
$$~~~~~~~~=\frac{1}{n}\int_{S^{n-1}}\left({\phi}
\left(\frac{\rho(K_{1},u)}{\rho(Q,u)}\right)+{\phi}
\left(\frac{\rho(L_{1},u)}{\rho(Q,u)}\right)\right)\rho(Q,u)
\rho(K_{2},u)\cdots\rho(K_{n},u)dS(u)$$
$$=\widetilde{V}_{\phi}(Q,K_{1},\cdots,K_{n})+\widetilde{V}_{\phi}
(Q,L_{1},K_{2},\cdots,
K_{n}).~~~~~~~~~~~~~~~~~~~~~~~~~~~~~~~~~\eqno(6.2)$$ Putting
$Q=K_{1}\widehat{+}_{{\phi}}L_{1}$ in (6.2), (6.2) changes (6.1).
~~~~~~~~~~~~~~~~~~~~~~~~~~~~~~~~~~~~~~~~~~~~~~~~~~~~~~~~~~~~~~~~~~~~~~~~~~~~~~~~~~~~~~~~~~~$\Box$

{\bf Theorem 6.2}~ (Dual Orlicz-Brunn-Minkowski inequality for the
Orlicz harmonic addition) {\it If $L_{1},K_{1},\cdots,K_{n}\in
{\cal S}^{n}$ and ${\phi}\in {\cal C}$, then for $\varepsilon>0$
$$\phi(1)\geq\phi\left(\frac{\widetilde{V}(K_{1},\ldots,K_{n})}{\widetilde{V}(K_{1}\widehat{+}_{\phi}\varepsilon\cdot L_{1},K_{2},\cdots,
K_{n})}\right)+\varepsilon\cdot\phi\left(\frac{\widetilde{V}(L_{1},K_{2},\cdots,
K_{n})}{\widetilde{V}(K_{1}\widehat{+}_{\phi}\varepsilon\cdot
L_{1},K_{2},\cdots, K_{n})}\right),\eqno(6.3)$$ If $\phi$ is
strictly convex, equality holds if and only if $K_{1}$ and $L_{1}$
are dilates.}

{\bf Proof}~ From Theorem 5.1 and Lemma 6.1, we have
$$\phi(1)\widetilde{V}(K_{1}\widehat{+}_{{\phi}}\varepsilon\cdot L_{1},K_{2},\cdots,K_{n})~~~~~~~~~~~~~~~~~~~~~~~~~~~~~~~~~~~~~~~~~~~~~~~~~~~~~~~~~~~~~~~~~~~~~~~~~~~~~~~~~~~~~~$$
$$~~~~~~~~~=\widetilde{V}_{\phi}(K_{1}\widehat{+}_{\phi}\varepsilon\cdot L_{1},K_{1},\cdots,K_{n})+\varepsilon\cdot\widetilde{V}_{\phi}
(K_{1}\widehat{+}_{\phi}\varepsilon\cdot L_{1},L_{1},K_{2},\cdots,
K_{n})~~~~~~~~~~~~~~~~~~~~~~~~~~~~~~~~~~~~~~~~~~~~~~~~~~~~~$$
$$~~~~~~~~~~~~~\geq\widetilde{V}(K_{1}\widehat{+}_{{\phi}}\varepsilon\cdot L_{1},K_{2},\cdots,K_{n})
\left\{\phi\left(\frac{\widetilde{V}(K_{1},\ldots,K_{n})}{\widetilde{V}(K_{1}\widehat{+}_{\phi}\varepsilon\cdot
L_{1},K_{2},\cdots,
K_{n})}\right)+\phi\left(\frac{\widetilde{V}(L_{1},K_{2},\cdots,
K_{n})}{\widetilde{V}(K_{1}\widehat{+}_{\phi}\varepsilon\cdot
L_{1},K_{2},\cdots, K_{n})}\right)\right \}.~~~~~~~$$ From the
equality condition of Theorem 5.1, the equality in (6.3) holds if
and only if $K_{1}$, $L_{1}$ and $K_{1}\widehat{+}_{\phi}L_{1}$
are dilates if $\phi$ is strictly convex, it follows that if
$\phi$ is strictly convex, the equality in (6.3) holds if and only
if $K_{1}$ and $L_{1}$ are dilates.
~~~~~~~~~~~~~~~~~~~~~~~~~~~~~~~~~~~~~~~~~~~~~~~~~~~~~~~~~~~~~~~~~~~~~~~~~~~~~~~~~~~~~~~~~~~~~~~~~~~~~~~~~~~~~~~~~~$\Box$

{\bf Theorem 6.3}~ (Dual Aleksandrov-Fenchel type inequality) {\it
If $L_{1},K_{1},\cdots,K_{n}\in {\cal S}^{n}$, $0\leq i,j<n$,
$1<r\leq n$ and ${\phi}\in {\cal C}$, then
$$\phi(1)\geq\phi\left(\frac{\prod_{i=1}^{r}\widetilde{V}(K_{i},\ldots,K_{i},K_{r+1},\ldots,K_{n})^{\frac{1}{r}}}{\widetilde{V}(K_{1}\widehat{+}_{\phi}L_{1},K_{2},\cdots,
K_{n})}\right)+\phi\left(\frac{\prod_{j=2}^{r}M(r)\cdot\widetilde{V}(K_{j},\ldots,K_{j},K_{r+1},\ldots,K_{n})^{\frac{1}{r}}}{\widetilde{V}(K_{1}\widehat{+}_{\phi}L_{1},K_{2},\cdots,
K_{n})}\right),\eqno(6.4)$$ where
$M(r)=\widetilde{V}(L_{1},\ldots,L_{1},K_{r+1},\ldots,K_{n})^{\frac{1}{r}}.$
If $\phi$ is strictly convex, equality holds if and only if
$L_{1}, K_{1},\ldots,K_{r}$ are all dilations of each other.}

{\bf Proof}~ This follows immediately from Theorem 6.2 and the
dual Aleksandrov-Fenchel inequality.
~~~~~~~~~~~~~~~~~~~~~~~~~~~~~~~~~~~~~~~~~~~~~~~~~~~~~~~~~~~~~~~~~~~~~~~~~~~~~~~~~~~~~~~~~~~~~~~~~~~~~~~~~~~~~~~~~~~~~~~~~~~~~~~~~~~~~~~~~~~~~~~~~~~~~~~$\Box$

{\bf Corollary 6.4}~ ($L_p$-dual Brunn-Minkowski inequality) {\it
If $L_{1},K_{1},\cdots,K_{n}\in {\cal S}^{n}$, $0\leq i,j<n$,
$1<r\leq n$ and $p\geq 1$, then
$$\widetilde{V}(K_{1}\widehat{+}_{p}L_{1},K_{2},\cdots,
K_{n})^{-p}\geq\prod_{i=1}^{r}\widetilde{V}(K_{i},\ldots,K_{i},K_{r+1},\ldots,K_{n})^{\frac{-p}{r}}~~~~~~~~~~~~~~~~~~~~~~~~~~~~~~~~~~~~~~~~$$
$$~~~~~~~~~~~~~~~~~+\prod_{j=1}^{r}M(r)^{-p}\cdot\widetilde{V}(K_{j},\ldots,K_{j},K_{r+1},
\ldots,K_{n})^{\frac{-p}{r}},\eqno(6.5)$$ with equality if and
only if $L_{1}, K_{1},\ldots,K_{r}$ are all dilations of each
other.}

{\bf Proof}~ This follows immediately from (6.4) with
$\phi(t)=t^{-p}$ and $p\geq 1$.
~~~~~~~~~~~~~~~~~~~~~~~~~~~~~~~~~~~~~~~~~~~~~~~~~~~~~~~~$\Box$

{\bf Corollary 6.5}~ {\it If $K,L\in {\cal S}^{n}$, ${\phi}\in
{\cal C}$ and $0\leq i<n-1$, then
$$\phi(1)\geq\phi\left(\left(\frac{\widetilde{W}_{i}(K)}{\widetilde{W}_{i}(K\widehat{+}_{\phi}L)}\right)^{\frac{1}{n-i}}\right)+\phi\left(\left(\frac{\widetilde{W}_{i}(L)}
{\widetilde{W}_{i}(K\widehat{+}_{\phi}L)}\right)^{\frac{1}{n-i}}\right).\eqno(6.6)$$
If $\phi$ is strictly convex, equality holds if and only if $K$
and $L$ are dilates.}

{\bf Proof}~ This follows immediately from Theorem 6.3 with
$r=n-i$, $K_{2}=\cdots=K_{n-i}=K\widehat{+}_{\phi}L$
$K_{n-i+1}=\cdots=K_{n}=B$.
~~~~~~~~~~~~~~~~~~~~~~~~~~~~~~~~~~~~~~~~~~~~~~~~~~~~~~~~~~~~~~~~~~~~~~~~~~~~~~~~~~~~~~~~~~~~~~~~~~~~~~~~$\Box$

The following inequality follows immediately from (6.6) with
$\phi(t)=t^{-p}$ and $p\geq 1$. If $K,L\in {\cal S}^{n}$, $0\leq
i<n$ and $p\geq 1$, then
$$\widetilde{W}_{i}(K\widehat{+}_{p}L)^{-p/(n-i)}\geq\widetilde{W}_{i}(K)^{-p/(n-i)}+\widetilde{W}_{i}(L)^{-p/(n-i)},$$
with equality if and only if $K$ and $L$ are dilates.

{\bf Corollary 6.6}~ {\it If $L_{1},K_{1},\cdots,K_{n}\in {\cal
S}^{n}$ and ${\phi}\in {\cal C}$, then
$$\phi(1)\geq\phi\left(\left(\frac{V(K_{1})\cdots
V(K_{n})}{\widetilde{V}(K_{1}\widehat{+}_{\phi}L_{1},K_{2},\cdots,
K_{n})^{n}}\right)^{\frac{1}{n}}\right)+\phi\left(\left(\frac{V(L_{1})V(K_{2})\cdots
V(K_{n})}{\widetilde{V}(K_{1}\widehat{+}_{\phi}L_{1},K_{2},\cdots,
K_{n})^{n}}\right)^{\frac{1}{n}}\right).\eqno(6.7)$$ If $\phi$ is
strictly convex, equality holds if and only if $L_{1},
K_{1},\ldots,K_{n}$ are all dilations of each other.}

{\bf Proof}~ This follows immediately from Theorem 6.3 with $r=n$.
~~~~~~~~~~~~~~~~~~~~~~~~~~~~~~~~~~~~~~~~~~~~~~~~~~~~~~~~$\Box$

{\bf Corollary 6.7}~ {\it If $L_{1},K_{1},\cdots,K_{n}\in {\cal
S}^{n}$ and $p\geq 1$, then
$$\widetilde{V}(K_{1}\widehat{+}_{p
}L_{1},K_{2}\cdots K_{n})^{-p}\geq(V(K_{1})\cdots
V(K_{n}))^{\frac{-p}{n}}+(V(L_{1})V(K_{2})\cdots V(
K_{n}))^{\frac{-p}{n}},\eqno(6.8)$$ with equality if and only if
$L_{1}, K_{1},\ldots,K_{n}$ are all dilations of each other.}

{\bf Proof}~ This follows immediately from (6.7) with
$\phi(t)=t^{-p}$ and $p\geq 1$.
~~~~~~~~~~~~~~~~~~~~~~~~~~~~~~~~~~~~~~~~~~~~~~~~~~~~~~~~$\Box$

Putting $K_{2}=\cdots=K_{n}=K_{1}\widehat{+}_{p}L_{1}$ in (6.8),
(6.8) becomes Lutwak's $L_{p}$-dual Brunn-Minkowski inequality
$$V(K\widehat{+}_{p}L)^{-p/n}\geq V(K)^{-p/n}+V(L)^{-p/n},$$
with equality if and only if $K$ and $L$ are dilates.

{\bf Corollary 6.8}~ {\it If $L_{1},K_{1},\cdots,K_{n}\in {\cal
S}^{n}$, $1\leq r\leq n$ and ${\phi}\in {\cal C}$, then
$$\widetilde{V}_{\phi}(L_{1},K_{1},K_{2},\cdots,
K_{n})\geq\widetilde{V}(L_{1},K_{2},\cdots, K_{n})
\phi\left(\frac{\prod_{i=1}^{r}\widetilde{V}(K_{i}\ldots,K_{i},K_{r+1},\ldots,K_{n})^{\frac{1}{r}}}{\widetilde{V}(L_{1},
K_{2}\ldots,K_{n})}\right).\eqno(6.8)$$ If $\phi$ is strictly
convex, equality holds if and only if $L_{1}, K_{1},\ldots,K_{r}$
are all dilations of each other.}

{\bf Proof}~ Let
$$K_{\varepsilon}=L_{1}\widehat{+}_{\phi}\varepsilon\cdot K_{1}.$$
From (4.8), dual Orlicz-Brunn-Minkowski inequality (6.3) and dual
Aleksandrov-Fenchel inequality, we obtain
$$\frac{1}{\phi'_{+}(1)}\cdot\widetilde{V}_{\phi}(L_{1},K_{1},\ldots,K_{n})
=\frac{d}{d\varepsilon}\bigg|_{\varepsilon=0^{+}}\widetilde{V}(K_{\varepsilon},K_{2},\cdots,K_{n})~~~~~~~~~~~~~~~~~~~~~~~~~~~~~~~~~~~~~~~~~~~~~~~~~~~~~~~~~~~$$
$$=\lim_{\varepsilon\rightarrow 0^{+}}\frac{\widetilde{V}(K_{\varepsilon},K_{2},\cdots,K_{n})-\widetilde{V}(L_{1},K_{2},\cdots,K_{n})}{\varepsilon}~~~~~~~~~~~~~~~~~~~~~~~~~~~~~~~~~~~~~~~~~~~~~~~~~~~~~~~$$
$$~~~~=\lim_{\varepsilon\rightarrow 0^{+}}\frac{\displaystyle
1-\frac{\widetilde{V}(L_{1},K_{2},\cdots,K_{n})}{\widetilde{V}(K_{\varepsilon},K_{2},\cdots,K_{n})}}
{\displaystyle
\phi(1)-\phi\left(\frac{\widetilde{V}(L_{1},K_{2},\cdots,K_{n})}{\widetilde{V}(K_{\varepsilon},K_{2},\cdots,K_{n})}\right)}\cdot\frac{\displaystyle
\phi(1)-\phi\left(\frac{\widetilde{V}(L_{1},K_{2},\cdots,K_{n})}{\widetilde{V}(K_{\varepsilon},K_{2},\cdots,K_{n})}\right)}{\varepsilon}\cdot\widetilde{V}(K_{\varepsilon},K_{2},\cdots,K_{n})$$
$$=\lim_{t\rightarrow 0^{+}}\frac{1-t}
{\displaystyle
\phi(1)-\left(t\right)}\cdot\lim_{\varepsilon\rightarrow
0^{+}}\frac{\displaystyle
\phi(1)-\phi\left(\frac{\widetilde{V}(L_{1},K_{2},\cdots,K_{n})}{\widetilde{V}(K_{\varepsilon},K_{2},\cdots,K_{n})}\right)}{\varepsilon}
\cdot\lim_{\varepsilon\rightarrow
0^{+}}\widetilde{V}(K_{\varepsilon},K_{2},\cdots,K_{n})~~~~~~~~~~$$
$$\geq\frac{1}{\phi'_{+}(1)}\cdot\lim_{\varepsilon\rightarrow
0^{+}}\phi\left(\frac{\widetilde{V}(K_{1},K_{2},\cdots,K_{n})}{\widetilde{V}(K_{\varepsilon},K_{2},\cdots,K_{n})}\right)
\cdot\widetilde{V}(L_{1},K_{2},\cdots,K_{n})~~~~~~~~~~~~~~~~~~~~~~~~~~~~~~~~~~~~~~~~$$
$$=\frac{1}{\phi'_{+}(1)}\cdot\phi\left(\frac{\widetilde{V}(K_{1},K_{2},\cdots,K_{n})}{\widetilde{V}(L_{1},K_{2},\cdots,K_{n})}\right)
\cdot\widetilde{V}(L_{1},K_{2},\cdots,K_{n})~~~~~~~~~~~~~~~~~~~~~~~~~~~~~~~~~~~~~~~~~~~~~~~~$$
$$\geq\frac{1}{\phi'_{+}(1)}\cdot\phi\left(\frac{\prod_{i=1}^{r}\widetilde{V}(K_{i}\ldots,K_{i},K_{r+1},\ldots,K_{n})^{\frac{1}{r}}}{\widetilde{V}(L_{1},
K_{2}\ldots,K_{n})}\right)\cdot\widetilde{V}(L_{1},K_{2},\cdots,K_{n}).~~~~~~~~~~~~~~~~~~~\eqno(6.9)$$
From (6.9), inequality (5.1) easy follows. From the equality
conditions of the dual Orlicz-Brunn-Minkowski inequality (6.3) and
dual Aleksandrov-Fenchel inequality, it follows that if $\phi$ is
strictly convex, the equality in (6.9) holds if and only if
$L_{1}, K_{1},\ldots,K_{n}$ are all dilations of each other.

This proof is
complete.~~~~~~~~~~~~~~~~~~~~~~~~~~~~~~~~~~~~~~~~~~~~~~~~~~~~~~~~~~~~~~~~~~~~~~~~~~~~~~~~~~~~~~~~~~~~~~~~~~~~~~~~~~~~~~~~$\Box$

\vskip 12pt

{\bf Availability of data and material}

All data generated or analysed during this study are included in
this published article.

{\bf Competing interests}

The author declare that he has no competing interests.

{\bf Funding}

The author's research is supported by the Natural Science
Foundation of China (11371334, 10971205).

{\bf Authors's contributions}

C-JZ contributed to the main results. The author read and approved the final manuscript.

{\bf Acknowledgements} The first author expresses his gratitude to
Professors G. Leng and W. Li for their valuable helps.

\end{document}